 \documentclass[final,3p,times,twocolumn]{elsarticle}
\usepackage{latexsym}
\usepackage{amsmath}%
\usepackage{amssymb}%
\usepackage{amsthm}
\usepackage{bm}
\usepackage{graphicx}
\usepackage{subfigure}
\graphicspath{{./},{./CMAME1/}}
\usepackage[english]{babel}
\usepackage[latin1]{inputenc}
\usepackage[T1]{fontenc}
\usepackage{enumerate} % to control the display of the enumeration counter
\theoremstyle{plain}

\newtheorem{lemma}{Lemma}

\newtheorem{proposition}{Proposition}
\newtheorem{remark}{Remark}

\numberwithin{equation}{section} % counter !!!!!!!
\newcommand\beq{\begin{equation}}
\newcommand\eeq{\end{equation}}
\renewcommand{\emph}{\textbf}
\newcommand{\brk}[1]{\left(#1\right)}          % \brk{.}     => (.)
\newcommand{\E}[1]{E\brk{#1}}
\newcommand{\Var}[1]{V\brk{#1}}
\newcommand{\D}{\mathcal{D}}

\newcommand{\xx}{\boldsymbol{x}}

\renewcommand{\to}{\rightarrow}

\newcommand{\N}{\mathbb{N}}
\newcommand{\R}{\mathbb{R}}

\newcommand{\mycomment}[1]{ }
\journal{\sep}

\begin{document}

%\markboth{SB}%S\'ebastien Boyaval}{Modelling and simulating}
\begin{frontmatter}

\title{A fast Monte-Carlo method with a Reduced Basis of Control Variates
% to Uncertainty Quantification
applied to Uncertainty Propagation and Bayesian Estimation} %Parameter Estimation}
%placed here. General acknowledgments should be placed at the end of the article.} % SPRINGER
%\subtitle{..} % SPRINGER
%\titlerunning{Reduced-Basis : Uncertain coefficients in elliptic PDEs}% SPRINGER

\author{S\'ebastien Boyaval
% \footnote{
% Laboratoire d'hydraulique Saint-Venant,\\ Universit\'e Paris-Est
% (Ecole des Ponts ParisTech) EDF R\&D,\\ 6 quai Watier, 78401 Chatou Cedex, France
% \\[2mm]       
% sebastien.boyaval@enpc.fr
% \\[2mm]
% Present address: EPFL
%}
}
%\authorrunning{Short form of author list} % SPRINGER if too long for running head

\address{
Universit\'e Paris Est, Laboratoire d'hydraulique Saint-Venant (EDF R\&D -- Ecole des Ponts ParisTech -- CETMEF), 
78401 Chatou Cedex, France ; and 
INRIA, MICMAC team--project, Domaine de Voluceau - Rocquencourt,
B.P. 105 - 78153 Le Chesnay, France.
E-mail: sebastien.boyaval@saint-venant.enpc.fr
}
%\institute{SB \at UPE 
%           \email{boyaval@}           %  \\
%           \emph{Present address:} of F. Author  %  if needed
%} % SPRINGER

% \begin{history}
% \received{}
% \revised{}
% \end{history}
%\date{Received: date / Accepted: date} % SPRINGER
%\date{\today}

%\maketitle

\begin{abstract}
The Reduced-Basis Control-Variate Monte-Carlo method was introduced recently in [S. Boyaval and T. Leli{\`e}vre, CMS, 8 2010] as an improved Monte-Carlo method, for the fast estimation of many parametrized expected values at many parameter values. We provide here a more complete analysis of the method including precise error estimates and convergence results. We also numerically demonstrate that it can be useful to some parametrized frameworks in Uncertainty Quantification, in particular
% It is based on a variance reduction technique so as to accelerate the reiterated computation of expectations.
% The principle is to build a Reduced Basis (RB) of control variates and next compute faster a parametrized expected value at many values of the parameter.
%We consider two kinds of applications in Uncertainty Quantification:
(i) the case where the parametrized expectation is a scalar output of 
the solution to a Partial Differential Equation (PDE) with stochastic coefficients (an Uncertainty Propagation problem), and
(ii) the case where the parametrized expectation is the Bayesian estimator of a scalar output in a similar PDE context.
% Numerous parametrized expectations at many parameter values are indeed often useful in Uncertainty Quantification frameworks (Uncertainty Propagation and Bayesian Estimation) 
% so as to: generate response surfaces, calibrate epistemic parameters in the presence of aleatoric parameters etc\ldots
Moreover, in each case, a PDE has to be solved many times for many values of its coefficients. This is costly and we also use a reduced basis of PDE solutions like in [S. Boyaval, C. Le Bris, Nguyen C., Y. Maday and T. Patera, CMAME, 198 2009]. This is the first combination of various Reduced-Basis ideas to our knowledge, here with a view to reducing as much as possible the computational cost of a simple approach to Uncertainty Quantification.
% numerical approximations using a Monte-Carlo/Finite-Element method. 
% \keywords{
% Parametrized Partial Differential Equations (PDEs) with random coefficients \and
% Reduced bases \and 
% Monte-Carlo \and 
% Variance reduction \and
% Control variates \and 
% Importance sampling \and
% Bayesian parameter estimation}
\end{abstract}

%\subjclass{} %


% \keywords{Finite element method, convergence analysis, existence of weak solutions.}%
% \ccode{AMS Subject Classification: 35Q30, 65M12, 65M60, 76A10, 76M10, 82D60}
 
\end{frontmatter}

%\tableofcontents

\section{Introduction}
\label{sec:intro}

The Reduced-Basis (RB) control-variate Monte-Carlo (MC) method was recently introduced in~\cite{boyaval-lelievre-2010} to compute fast many expectations of scalar outputs of the solutions to parametrized ordinary Stochastic Differential Equations (SDEs) at many parameter values. But as a simple, generic MC method with reduced variance, the RB control-variate MC method can also be useful in other parametric contexts and the main goal of this article is to show that it can be useful to Uncertainty Quantification (UQ) too, possibly in combination with the standard RB method in a PDE context.

There is a huge literature on the UQ subject. Indeed, to be actually predictive in real-life situations~\cite{liu-bayarri-berger-paulo-sacks-2008}, most numerical models require (i) to calibrate as much as possible the parameters and (ii) to quantify the remaining uncertainties propagated by the model. Besides, the latter two steps are complementary in an iterative procedure to improve numerical models using datas from experiments: quantifying the variations of outputs generated by input parameters allows one to calibrate the input uncertainties with data and in turn reduces the epistemic uncertainty in outputs despite irreducible aleatoric uncertainty.
% aleatoric = statistical / can be estimated from long time seuqnces as frequences, while epistemic = systematic / depend on the model
% \footnote{
% Statistical aleatoric uncertainty seems unavoidable whatever the scale described by a model
% since most elementary physical principles like quantum mechanics are probabilistic in essence.
% }.
Various numerical techniques have been developped to quantify uncertainties and have sometimes been used for years~\cite{ghanem-spanos-1991,bergman-et-al-1997}. % a review of some methods in the introduction of the monograph
But there are still a number of challenges~\cite{debusschere-najm-pebay-knio-ghanem-lemaitre-2004,burkardt-gunzburger-webster-2007,doostan-ghanem-redhorse-2007}. 

For PDEs in particular, the coefficients are typical sources of uncertainties.
% There are basically two kinds of random models for uncertainties:
% (i) the random input (here, stochastic coefficients in a PDE) is endowed 
% with a probability distribution that presumably belongs to 
% some parametric family where the family parameters -- often termed ``hyper-parameters'' -- need calibrating
% (see e.g.~\cite{casella-berger-1990,xiu-karniadakis-2003}), or 
% (ii) the random input is estimated by means of general non-parametric distributions,
% defined either directly at the computer level using any distribution or first abstractely using Dirichlet process
% (see e.g.~\cite{tsybakov-2009}).
One common modelling of these uncertainties endows the coefficients with a probability distribution that presumably belongs to some parametric family and the PDEs solutions inherit the random nature of the uncertainty sources. A Bayesian approach is often favoured to calibrate the parameters in the probability law using observations of the reality~\cite{wang-zabaras-2005,marzouk-najm-2009}. But the accurate numerical simulation of the PDEs solutions as a function of parametrized uncertain coefficients is a computational challenge due to its complexity, and even more so is the numerical optimization of the parameters in uncertain models. That is why new/improved techniques are still being investigated~\cite{ghanem-doostan-2006,soize-2008}. % soize-ghanem-2004
Our goal in this work is to develop a practically useful numerical approach that bases on the simple MC method to simulate the probability law of the uncertain coefficients. We 
suggest to use the RB control-variate MC method in some UQ frameworks to improve the computational cost of the naive MC method, in particular in contexts where some coefficients in a PDE are uncertain and other are controlled.

There exist various numerical approaches to UQ. The computational cost of MC methods is certainly not optimal when the random PDE solution is regular, see e.g.~\cite{cohen-devore-schwab-2010}. But we focus here on MC methods because they are (a) very robust, that is useful when the regularity of the solution with respect to the random variables degrades, and (b) very easy to implement (they are {\it non-intrusive} in the sense that they can use a PDE numerical solver as a black-box, with the values of the PDE coefficients as input and that of the discrete solution as output). Besides, note that even when the random PDE solution is very regular with respect to the random variables, it is not yet obvious how algorithms can take optimal profit of the regularity of random PDE solutions and remain practically efficient as the dimension of the (parametric) probability space increases, see e.g.~\cite{cohen-devore-schwab-2011}. So, focusing on a MC approach, our numerical challenge here is basically two-sided: (i) on the probabilistic side, one should sample fast the statistics of the random PDE solution (or of some random output that is the quantity of interest), and (ii) on the deterministic side, one has to compute fast the solution to a PDE for many realizations of the random coefficients. It was proposed in~\cite{boyaval-lebris-maday-nguyen-patera-2009} to use the RB method in order to reduce the numerical complexity of (ii), but this does not fully answer the numerical challenge. In particular, although the RB method can improve naive MC approaches at no-cost (since the selection of the reduced basis for the PDE solutions at various coefficients values can be trained on the same large sample of coefficients values that is necessary to the MC sampling), the resulting MC approach might still be very costly, maybe prohibitively, due to the large number of realizations that is necessary to accurately sample the statistics of the PDE solution (side (i) of our challenge above). In this work, we thus tackle the question how to reduce the numerical complexity of (i). We have in mind the particular but useful case where one is interested in the expected value of a random scalar output of the random PDE solution as a function of a (deterministic) {\it control parameter}, typically another (deterministic) coefficient in the UQ problem which is ``under control''. (Think of the construction of response surfaces for a mean value as a function of control parameters.) A similar parametric context occurs in Bayesian estimation, sometimes by additionally varying the {\it hyper parameters} or the observations. In any case, our goal is to reduce the computational cost of a parametrized (scalar) MC estimation when the latter has to be done many times for many values of a parameter, and we illustrate it with examples meaningful in a UQ context.

To accelerate the convergence of MC methods as numerical quadratures for the expectation of a random variable, one idea is to tune the sampling for a given family of random variables like in the quasi-Monte-Carlo (qMC) methods~\cite{sloan-wozniakowski-1998,venkiteswaran-junk-05-1,kuo-schwab-sloan-2011}. Another common idea is to sample another random variable with same mean but with a smaller variance. % than the one whose first moment we are interested in.
Reducing the variance allows one to take smaller MC samples of realizations and yet get MC estimations with confidence intervals of similar (asymptotic) probability. Many techniques have been designed in order to reduce the variance in various contexts~\cite{hammersley-handscomb-1965,jourdain-2009}. Our RB control-variate MC method bases on the so-called control-variate technique. It has a specific scope of application in parametric contexts. But it suits very well to some computational problems in mathematical finance and molecular dynamics as shown in~\cite{boyaval-lelievre-2010}, and can be useful in UQ as we are going to see. 
 
The paper is organized as follows.
In Section~\ref{sec1}, we recall the RB control-variate technique as a general variance reduction tool
for the MC approximation of a parametrized expected value at many values of the parameter. The presentation is a bit different to that in~\cite{boyaval-lelievre-2010}, which was more focused on SDEs. Moreover, we also give new error estimates and convergence results.
% The principle of our variance reduction technique is based on the same paradigm as the standard RB method, hence its name. 
% It can quickly find a good control variate for a parametrized expected value at any value of the parameter 
% (in fact, a good linear combination of precomputed control variates that define a linear space) 
% at the price of precomputations (often made in a pre-processing stage termed ``offline'') 
% and a few additional computations (done ``online'', according to the terminology of algorithm with real-time constraints).
% Since precomputations are needed to identify a good basis for control variates in a trial set of expected values, one needs of course to identify a regime of efficiency.
% As illustrative examples for application in UQ frameworks, we consider two particular cases where one is interested in the variations with respect to a parameter of an expected value 
% that is some output of a solution to a PDE with stochastic coefficients. 
% First, a case where the parameter entering the output expected value is low-dimensional, namely a few other -- deterministic -- coefficients in the PDE 
% (think of sensitivity analyses or of the construction of response surfaces in UQ as a function of a control parameter).
% Second, a case where it is high-dimensional, namely coefficients that directly enter the output expected value
% (think of Bayesian estimation as a function of a control parameter and a hyper parameter that has to be adjusted).
In Section~\ref{sec2}, the RB control-variate MC method is applied to compute the mean of a random scalar output
in a model PDE with stochastic coefficients (the input uncertainty) at many values of a control parameter.
In Section~\ref{sec3}, it is applied to Bayes estimation, first for a toy model where various parametric contexts are easily discussed, then for the same random PDE as in section~\ref{sec2}. 

We also note that this work does not only improve on the RB approach to UQ~\cite{boyaval-lebris-maday-nguyen-patera-2009} but also on an RB approach to Bayesian estimation proposed in~\cite{nguyen-rozza-huynh-patera-2009} with a deterministic quadrature formula to evaluate integrals. 
% thus assuming enough regularity in the random variations to be computationally more efficient than a simple MC
For both applications, to our knowledge, our work is the first attempt at optimally approximating the solution with a simple MC/FE method by combining RB ideas of two kinds, stochastic and deterministic ones~\cite{boyaval-lebris-lelievre-maday-nguyen-patera-2010}. Note that for convenience of the reader non-expert in RB methods, the standard RB method~\cite{maday-2006,patera-rozza-2007} is briefly recalled in Section~\ref{sec:certifiedRBmethod}.

% \begin{remark}[On the choice of the model problems]
% In the numerical applications, we propagate a one-dimensional stochastic coefficient throug a scalar PDE parametrized by another one-dimensional coefficient,
% and we compute the Bayesian Minimum-Mean-Square-Error (MMSE) estimator of a one-dimensional parameter from a sample of realizations of a stochastic model on the other hand.
% In both cases, the outputs are parametrized by other coefficients ``under control'' (not necessarily one-dimensional). Thus they define simple enough many-query contexts
% to show the applicability of our ideas in some UQ frameworks and to motivate further investigations
% (in particular, a more sensible -- but more involved ! -- implementation for high-dimensional problems).
% Of course we are aware that the MC method is certainly not optimal for sampling one-dimensional integrals,
% and that more work is necessary to demonstrate the full applicability of our ideas in the most realistic UQ problems.
% \end{remark}

%\section{The Reduced-Basis Control-Variate Monte-Carlo Method}
\section{The RB Control-Variate MC Method}
\label{sec1}

The RB control-variate technique is a generic variance reduction tool for the MC approximation of a parametrized expected value at many values of the parameter. In this section we recall the technique for the expectation $E(Z^{\lambda})$ of a generic square-integrable random variable $Z^{\lambda}\in L^2_P$ parametrized by $\lambda$. The principle for the reduction of computations is based on the same paradigm as the standard RB method and allows one to accelerate the MC computations of many $E(Z^{\lambda})$ at many values of $\lambda$. Our presentation is slightly different than the initial one in~\cite{boyaval-lelievre-2010} and gives new elements of analysis (error estimates and convergence results).
% There we showed it was efficient for
% (i) the calibration of a parametric volatility in Black-Scholes equation, and
% (ii) the multiscale computation of viscoelastic flows (one instance of molecular dynamics 
% coupled with a macroscopic mechanical evolution equation).
% We perform a deeper analysis of the method than the one in~\cite{boyaval-lelievre-2010}.

\subsection{Principles of RB control-variate MC method}

Let $P$ be a probability measure such that $Z^{\lambda}$ is a random variable in $L^2_P$ for all parameter values $\lambda$ in a given fixed range $\Lambda$. Assume one has an algorithm to simulate the law of $Z^{\lambda}$ whatever $\lambda\in\Lambda$. Then, at any $\lambda\in\Lambda$, one can define MC estimators $E_M(Z^{\lambda})$ that provide useful approximations of $E(Z^{\lambda})$, by virtue of the strong law of large numbers
\beq
\label{eq:MC0}
E_M(Z^{\lambda}) := \frac1M \sum_{m=1}^M Z^{\lambda}_m \xrightarrow[M\to\infty]{P-a.s.} E(Z^{\lambda}) \,,
\eeq
provided the number $M$ of independent identically distributed (i.i.d.) random variables $Z^{\lambda}_m \sim Z^{\lambda}$, $m=1\ldots M$, is sufficiently large. Here, the idea is: if $E(Z^{\lambda_i^I}):=\int Z^{\lambda_i^I} dP$ is already known with a good precision for $I$ parameter values $\lambda_i^I$, $i=1\ldots I$, ($I\in\N_{>0}$) and if the law of $Z^{\lambda}$ depends smoothly on $\lambda$ then, given $I$ well-chosen real values $\alpha_i^\lambda$, $i=1\ldots I$, the standard MC estimator $E_M(Z^{\lambda})$ could be efficiently replaced by a MC estimator for $E(Z^{\lambda}-\sum_{i=1}^I\alpha_i^\lambda Z^{\lambda_i^I})+\sum_{i=1}^I\alpha_i^\lambda E(Z^{\lambda_i^I})$ that is as accurate and uses much less than $M$ copies of $Z^{\lambda}$.

In other words, if the random variable % $Y^{\lambda}$
\beq
\label{eq:controlvar}
\hat Y^\lambda = \sum_{i=1}^I \alpha_i^\lambda (Z^{\lambda_i^I}-E(Z^{\lambda_i^I}))
\eeq
is correlated with $Z^{\lambda}$ such that the control of $Z^{\lambda}$ by $\hat Y^{\lambda}$ reduces the variance, that is if
$$ V(Z^{\lambda}) := \int |Z^\lambda-E(Z^{\lambda})|^2 dP \ge V(Z^{\lambda}-\hat Y^{\lambda}) \,,$$
then the confidence intervals with asymptotic probability $\mathop{\rm erf}\left(\frac{a}{\sqrt{2}}\right)$ ($a>0$) for MC estimations
\begin{equation}
%\begin{multline}
\label{eq:MC}
E_M(Z^{\lambda}-\hat Y^{\lambda}) := \sum_{m=1}^M \frac{Z^{\lambda}_m-\hat Y^{\lambda}_m}M
%\\ 
\xrightarrow[M\to\infty]{P-a.s.} %E(Z^{\lambda}-Y^{\lambda}) \equiv 
E(Z^{\lambda})
%\end{multline}
\end{equation}
that are in the Central Limit Theorem (CLT)
\beq
%\begin{multline}
\label{eq:CLT}
P\left( 
 % \left| E_M((Z^{\lambda}-Y^{\lambda})_m) - E(Z^{\lambda}-Y^{\lambda}) \right| \le a\sqrt{\frac{V(Z^{\lambda}-Y^{\lambda})}{M}} 
 \frac{\left| E_M(Z^{\lambda}-\hat Y^{\lambda}) - E(Z^{\lambda}) \right|}{\sqrt{V(Z^{\lambda}-\hat Y^{\lambda})/M}} \le a 
\right)
%\\ 
\underset{M\to\infty}{\longrightarrow} %\int_{-a}^a \frac{e^{-x^2/2}}{\sqrt{2\pi}} dx \equiv 
\mathop{\rm erf}\left(\frac{a}{\sqrt{2}}\right)
%\end{multline}
\eeq
converge faster with respect to the number $M$ of realizations than the confidence intervals for~\eqref{eq:MC0}. The unbiased estimator~\eqref{eq:MC} ($E(\hat Y^{\lambda})=0$) is thus a better candidate than~\eqref{eq:MC0} for a fast MC method.

At a given $\lambda\in\Lambda$, we define $\alpha_i^\lambda$, $i=1\ldots I$, in~\eqref{eq:controlvar} to obtain the optimal variance reduction minimizing $V(Z^{\lambda}-\hat Y^{\lambda})$. Equivalently, the control variates of the form~\eqref{eq:controlvar} are thus defined as best approximations of the ideal control variate $Y^\lambda:=Z^{\lambda}-E(Z^{\lambda})\in L^2_P $ that reduces the variance to zero%, viz. solutions to least-squares problems
\begin{equation}
\label{eq:variance_minimization}
%\underset{(\alpha_i^\lambda)\in\R^I}{\mathop{\rm inf}}\:\Var{Z^\lambda-Y^\lambda} \equiv
\underset{(\alpha_i^\lambda)\in\R^I}{\mathop{\rm inf}}
%E\left( \left| Z^{\lambda}-E(Z^{\lambda}) -\sum_{i=1}^I\alpha_i^\lambda(Z^{\lambda_i^I}-E(Z^{\lambda_i^I}))\right|^2 \right) \,.
E\left( \left| Y^{\lambda}-\sum_{i=1}^I\alpha_i^\lambda Y^{\lambda_i^I}\right|^2 \right) \,.
\end{equation}
The $\alpha_i^\lambda$, $i=1\ldots I$, thus solve a least-squares problem with normal equations for $i=1\ldots I$
\begin{equation}
\label{eq:RBnormal}
\sum_{j=1}^I C\left(Z^{\lambda_i^I},Z^{\lambda_j^I}\right) \alpha_j^\lambda = C\left(Z^{\lambda_i^I},Z^{\lambda}\right) \,,
\end{equation}
where $C(Y,Z)=E((Y-E(Y))(Z-E(Z)))$ denotes the covariance between $Y,Z\in L^2_P$. %two square-integrable random variables $Y$ and $Z$,
And they are unique as long as the %positive 
matrix $\mathbf{C}$ with entries $\mathbf{C}_{i,j}=C\left(Z^{\lambda_i^I},Z^{\lambda_j^I}\right)$ remains definite. % (thus, invertible).

In general, the computation of a good control variate is difficult (the ideal one $Y^{\lambda}$ requires the result $E(Z^{\lambda})$ itself). But we proved in~\cite{boyaval-lelievre-2010} that control variates of the form~\eqref{eq:controlvar} actually make sense in some contexts, a result inspired by~\cite{maday-patera-turinici-2002a,maday-patera-turinici-2002b}. Let us denote by $L^2_{P,0}$ the Hilbert subspace of random variables $Y\in L^2_P$ that are centered $E(Y)=0$.
\begin{proposition}
\label{prop:apriori_rb}
Consider a set of random variables % parametrized by $\lambda \in \Lambda$
\begin{equation} \label{eq:Yassumption}
Y^\lambda = \sum_{j=1}^J g_j(\lambda) \, Y_j \,,\ \forall \lambda \in \Lambda \,,
\end{equation}
where $Y_j \in L^2_{P,0}$, $j=1\ldots J$, are uncorrelated % parameter-independent random variables  
and $(g_j)_{1\le j\le J}$ are $C^\infty_{\ge0}(\R)$ % positive %real analytic 
functions. % exist such that % \in C^\infty(\Lambda,\R_{>0})
If there exists a constant $C>0$ and an interval $\widetilde\Lambda$ such that, 
for all parameter ranges $\Lambda = [\lambda_{\min},\lambda_{\max}] \subset \R$,
there exists a $C^\infty$ diffeomorphism $\tau_\Lambda:\Lambda\to\widetilde\Lambda$ where %it holds
\begin{equation}\label{eq:gassumption}
\sup_{1\le j\le J} \sup_{\tilde{\lambda}\in\tau_{\Lambda}(\Lambda)}
 (g_j\circ\tau_\Lambda^{-1})^{(M)}(\tilde{\lambda}) \le M! C^M \,, 
%\text{ for all $M$-derivatives of $g_j\circ\tau_\Lambda^{-1}$,} 
\end{equation}
for all $M$-derivatives $(g_j\circ\tau_\Lambda^{-1})^{(M)}$ of $g_j\circ\tau_\Lambda^{-1}$,
then there exist constants $c_1,c_2>0$ independent of $\Lambda$ and $J$ such that, 
for all parameter ranges $\Lambda = [\lambda_{\min},\lambda_{\max}] \subset \R$, %the following bound holds 
for all $I\in\mathbb{N}$, $ I\ge I_0 >0 $,
% $$ N_0 = O\left(\ln\left(\frac{\lambda_{\max}}{\lambda_{\min}}\right)\right) \textit{ as } \to+\infty \,, $$
\begin{equation} \label{eq:rbcv} % \sup_{\lambda\in\Lambda} 
\inf_{Y \in \mathcal{Y}_I} V(Y^\lambda-Y) \le 
e^{-\frac{c_2(I-1)}{I_0-1}} \, 
%\left( \sum_{j=1}^J |g_j(\lambda)|^2\:\Var{Y_j}\right)
V(Y^\lambda) 
\,,\ {\forall\lambda\in\Lambda} \,,
\end{equation}
where $I_0 := 1+c_1\left(\tau_{\Lambda}(\lambda_{\max})-\tau_\Lambda(\lambda_{\min})\right) $
and $\mathcal{Y}_I=\mathop{\rm span}(Y^{\lambda_i^I},i=1,\ldots,I)$ 
%invokes $I$ explicit parameter values $ \lambda_i^I\in\Lambda\,,\ i=1,\ldots,I $ (see~\cite{boyaval-lelievre-2010}).
uses the $I$ explicit values
$\lambda_i^I = \tau_\Lambda^{-1}\left(\tau_\Lambda(\lambda_{\min}) + \frac{i-1}{I-1}\left(\tau_\Lambda(\lambda_{\max})-\tau_\Lambda(\lambda_{\min})\right)\right)$.
\end{proposition}
\noindent 
Prop.~\ref{prop:apriori_rb} tells us that using~\eqref{eq:controlvar} as a control variate makes sense because the variance in~\eqref{eq:rbcv} decays very fast with $I$ for all $\lambda\in\Lambda$, and whatever $\Lambda$. %(a Gram-Schmidt procedure allows one to write $Y^\lambda$ like in~\eqref{eq:Yassumption}, using uncorrelated random variables, at least on small ranges $\Lambda$ %where they do not vanish, for the coefficients to remain positive).

But the explicit construction of $\lambda_i^I$, $i=1\ldots I$, in Prop.~\ref{prop:apriori_rb}, %of good control variates 
for random variables $Z^{\lambda}$ that are analytic in a 1D parameter $\lambda$, does not generalize to higher-dimensions in an efficient way, a well-known manifestation of the {\it curse of dimensionality}
\footnote{
 A tensor-product of~\eqref{eq:rbcv} in dimension $d$ yields a rate $-I/d$.
}.
Yet, our numerical results (here and in~\cite{boyaval-lelievre-2010}) show our variance reduction principle can work well in contexts with higher-dimensional or less regular parametrization.
% even when the very restrictive real-analyticity assumption~\eqref{eq:gassumption} is not satisfied.
Then, given a parametric family $Z^{\lambda}$, $\lambda\in\Lambda$, how to choose parameter values $\lambda_i^I$, $i=1\ldots I$, that can in practice efficiently reduce the MC computations at all $\lambda\in\Lambda$ ?

Given $I\in\mathbb{N}$, let us define $\lambda_i^I$, $i=1\ldots I$, as the parameter values with optimal approximation properties in $L^\infty_{\Lambda}(L^2_{P,0})$: the linear space spanned by $Y^{\lambda_i^I}=Z^{\lambda_i^I}-E(Z^{\lambda_i^I})$ minimizes the variances over all $\lambda\in\Lambda$, or equivalently the $\lambda_i^I$, $i=1\ldots I$, minimize
% of $\{ Z^{\lambda}-E(Z^{\lambda}), \lambda \in \Lambda\}$
\begin{equation}
\label{eq:bestNvariates}
\bar d_I:=\underset{\{\lambda_1^I,\ldots,\lambda_I^I\}\in\Lambda^I}{\mathop{\rm inf}}\:\
%\left(
\underset{\lambda\in\Lambda}{\sup}\:\
\underset{(\alpha_i^\lambda)\in\R^I}{\mathop{\rm inf}}\:\
V(Z^\lambda-\hat Y^\lambda) \,.
%\right)
\end{equation}
The variates $Y^{\lambda_i^I}$, $i=1\ldots I$, define an optimal $I$-dimensional {\it reduced basis} in $L^2_{P,0}$ for the approximations $\hat Y^\lambda$ of the ideal controls $Y^{\lambda}$ at all $\lambda\in\Lambda$. The values $\bar d_I, I\in\N,$ in~\eqref{eq:bestNvariates} coincide with a well-known notion in approximation theory~\cite{pinkus-1985}:
the Kolmogorov widths of $\mathcal{Y}:=\{ Z^{\lambda}-E(Z^{\lambda}), \lambda\in\Lambda\}$ in $\mathcal{Y}\subset L^2_{P,0}$, and are upper bounds for the more usual Kolmogorov widths of $\mathcal{Y}$ in $L^2_{P,0}$ % i.e. the worst error commited by the best $I$-linear approximation in $L^2_{P,0}$
$$
\begin{array}{l}
d_I :=
\underset{\left\{Y_i^I\right\}\in [L^2_{P,0}]^I}{\inf}\:\
\underset{\lambda\in\Lambda}{\sup}\:\
\underset{Y\in \mathop{\rm span}\left(Y_1^I,\ldots,Y_I^I\right)}{\inf}\:\
\| Y^\lambda-Y \|_{L^2_{P}} \,.
\end{array}
$$

Explicitly computing Kolmogorov widths is a difficult problem solved only for a few simple sets $\mathcal{Y}$, let alone the computation of (non-necessarily unique) minimizers. But in a {\it many-query} parametric framework where the MC computations are queried many times for many values of the parameter, one can take advantage of a presumed smooth dependence of the MC computations on the parameter to identify a posteriori some approximations for the optimal parameter values $\lambda_i^I$, $i=1\ldots I$, with prior computations at only a few parameter values. Moreover, in practice, we shall be content with approximations of $\lambda_i^I$, $i=1\ldots I$ that are not very good, provided they produce sufficiently fast-decaying upper-bounds of $\bar d_I$ as $I\to\infty$.
% see e.g.~\cite{rudin-1952,pinkus-1985,buslaev-tikhomirov-1992,izaak-1994,
% kushpel-2000,magaril-osipenko-tikhomirov-2001,uskov-2002,kudryavtsev-2005}.
%\text{ if $I>0$ and } d_0(\Lambda) := \lim_{M_{\rm small}\to\infty} \lim_{M_{\rm large}\to\infty} \sigma_{0,M_{\rm small},M_{\rm large}}(\Lambda)
To our knowledge, this reduction paradigm has been applied for the first time in~\cite{boyaval-lelievre-2010} to minimize the variance of many parametrized MC estimations. % with control variates~\eqref{eq:controlvar}.
It is inspired by the (by now) standard RB method to minimize the discretization error in many-query Boundary Value Problems (BVP) for parametrized PDEs (see~\cite{maday-2006,patera-rozza-2007} e.g. and Section~\ref{sec:certifiedRBmethod}). %developped by Y.Maday, A.T.Patera and collaborators.

Let us define linear combinations like~\eqref{eq:controlvar} when $I=0$ by $0$. We next define a (non-unique) sequence of parameter values $\lambda_i\in\Lambda$, $i\in\mathbb{N}_{>0}$ by the sucessive iterations $I\in\mathbb{N}$ of a {\it greedy} algorithm
\beq
%\begin{multline}
\label{eq:greedy} %\lambda_1 \underset{\lambda\in\Lambda}{\mathop{\rm argsup}}Var(Z^\lambda) \\
\lambda_{I+1}\in\underset{\lambda\in\Lambda}{\mathop{\rm argsup}}V(Z^\lambda-\sum_{i=1}^IZ^{\lambda_i}) \,.
%\end{multline}
\eeq
For all $I\in\mathbb{N}$, the parameter values $\lambda_i$, $i=1\ldots I$, approximate $\lambda_i^I$, $i=1\ldots I$ in~\eqref{eq:bestNvariates} and yield the following upper-bound
\beq
%\begin{multline}
\label{eq:greedyres} %\lambda_1 \underset{\lambda\in\Lambda}{\mathop{\rm argsup}}Var(Z^\lambda) \\
\sigma_I:=\underset{\lambda\in\Lambda}{\mathop{\rm sup}}V(Z^\lambda-\sum_{i=1}^I\alpha_i^\lambda Z^{\lambda_i}) 
\ge \bar d_I \ge d_I \,.
%\end{multline}
\eeq
Using a finite discrete subset $\tilde\Lambda\subset\Lambda$ (possibly $\tilde\Lambda=\Lambda$ in case of finite cardinality ${\rm card}(\Lambda)<\infty$) and a MC estimator for the variance, like
\begin{multline}
\label{eq:varianceMCestimation}
V_M(Z^\lambda-\sum_{i=1}^IZ^{\lambda_i}) %-E(Z^{\lambda_i}))
% = \frac{\sum_{m=1}^M \left| (Z^\lambda_m-\sum_{i=1}^IZ^{\lambda_i}_m)-E_m(Z^\lambda-\sum_{i=1}^IZ^{\lambda_i}) \right|^2}
% {(M-1)\left(1-\sum_{m=1}^{M-1} \frac1{M(m+1)}\right)} \,,
= \frac{1}{M-1}
\sum_{m=1}^M \bigg| (Z^\lambda_m-\sum_{i=1}^I\alpha_i^\lambda Z^{\lambda_i}_m) 
\\ - \frac{m-1}{m} E_m(Z^\lambda-\sum_{i=1}^I\alpha_i^\lambda Z^{\lambda_i}) \bigg|^2 \,,
\end{multline}
we can finally define a computable sequence of parameter values $\tilde\lambda_i\in\Lambda$, $i\in\mathbb{N}_{>0}$ by sucessive iterations $I\in\mathbb{N}$ of a {\it weak greedy} algorithm at $\omega\in\Omega$
\beq
%\begin{multline}
\label{eq:weakgreedy} %\lambda_1 \underset{\lambda\in\Lambda}{\mathop{\rm argsup}}Var(Z^\lambda) \\
\tilde\lambda_{I+1}\in\underset{\lambda\in\tilde\Lambda}{\mathop{\rm argsup}}V_M(Z^\lambda-\sum_{i=1}^IZ^{\tilde\lambda_i})(\omega) \,.
%\end{multline}
\eeq

For all $I\in\mathbb{N}$, our computable approximations of the optimal parameter values $\lambda_i^I$, $i=1\ldots I$ used hereafter in practical applications will be $\tilde\lambda_i\equiv\tilde\lambda_i(\omega)$, $i=1\ldots I$, from a set like in~\eqref{eq:weakgreedy}. Note that the $\tilde\lambda_i$, $i\in\N_{>0}$, are constructed iteratively and thus do not depend on $I$ but only on $\tilde\lambda_1$ (and possibly on $\tilde\Lambda$), which is useful in practice: the adequate value of $I$ for a good variance reduction is typically not known in advanced and one simply stops the weak greedy algorithm at $I\in\mathbb{N}$ when, at a given $\omega\in\Omega$,
\beq
%\begin{multline}
\label{eq:weakgreedyres} %\lambda_1 \underset{\lambda\in\Lambda}{\mathop{\rm argsup}}Var(Z^\lambda) \\
\tilde\sigma_I:=\underset{\lambda\in\tilde\Lambda}{\mathop{\rm sup}}V_M(Z^\lambda-\sum_{i=1}^I\alpha_i^\lambda Z^{\tilde\lambda_i}) (\omega)\le TOL
%\end{multline}
\eeq
is smaller than a maximum tolerance $TOL>0$. 

Even if greedy algorithms might yield approximate minimizers of~\eqref{eq:bestNvariates} that are far suboptimal, numerical results show that they can nevertheless be useful to computational reductions using the RB control-variate MC method, just as they already proved useful to computational reductions with the standard RB method in numerous practical examples. Recent theoretical results~\cite{buffa-maday-patera-prudhomme-turinici-2009,binev-cohen-dahmen-devore-petrova-wojtaszczyk-2010} also support that viewpoint as concerns the standard RB method and can be straightforwardly adapted to our framework. Comparing directly $\sigma_I$ with $d_I$ at same $I\in\N$ will not, in general, give estimates better than $\sigma_I\le \frac{2^{I+1}}{\sqrt{3}}d_I$. And this is pessimistic as regards the convergence of the greedy algorithm insofar as it predicts variance reduction only for sets with extremely fast decaying Kolmogorov widths $d_I$. Yet, the two sequences $\sigma_I$ and $d_I$ typically have comparable decay rates as $I\to\infty$, and it holds for adequate $\eta,\beta,c,C>0$ given any $d_0,\alpha,a>0$:
\begin{itemize}
 \item[a)] $d_I \le d_0 I^{-\alpha}\,, \forall I \Rightarrow \sigma_{I} \le C d_0 I^{-\alpha}\,, \forall I$,
 \item[b)] $d_I \le d_0 e^{-a I^{-\alpha}}\,, \forall I \Rightarrow \sigma_I \le C d_0 e^{-c I^{-\beta}}\,, \forall I $,
 \item[c)] $d_I \le d_0 e^{-a I^{-\alpha}}\,, \forall I \Rightarrow \sigma_I \le C I \eta^I e^{-a I^{-\alpha}}\,, \forall I$,
\end{itemize}
where c) is sharper than b) if, and only if, $\alpha>1$ or $\alpha=1$ and $a>\ln2$. So, when the Kolmogorov widths $d_I$ decay fast (variance reduction is a priori possible), greedy algorithms can be useful if used with $I$ sufficiently many iterations
\footnote{
 Then computational reductions are possible, but only for sufficiently many queries in the parameter of course, as we shall see.
}.

The latter results can also be adapted to the weak greedy algorithm, like in~\cite{binev-cohen-dahmen-devore-petrova-wojtaszczyk-2010}. For all $\epsilon\in(0,1)$ and $I\in\mathbb{N}$, given any $\theta_I\in(0,1)$, we define the smallest number of realizations $M_{\rm test}^{\epsilon,I}$ such that
$$
P\left( |V_{M_{\rm test}^{\epsilon,I}}(Z^\lambda-\hat Y^{\lambda})-V(Z^\lambda-\hat Y^{\lambda})|\le \theta_Id_I \right) \ge 1-\epsilon
$$
holds for all $\lambda$ and control variates $\hat Y^{\lambda}$ of dimension $0\le i\le I$.
Of course, we do not know exactly $M_{\rm test}^{\epsilon,I}$ in practice, which could even not be finite if variations in $\lambda$ are not smooth enough. But if, for all $\epsilon\in(0,1)$, we assume that one can use more than $M_{\rm test}^{\epsilon,I}$ realizations in the computations of step $I$ of the weak greedy algorithm, then the following events occur with probability more than $1-\epsilon$, for adequate $\eta,\beta,c,C>0$ given any $d_0,\alpha,a>0$:
\begin{itemize}
 \item[a)] $d_I \le d_0 I^{-\alpha}\,, \forall I \Rightarrow \tilde\sigma_{I} \le C d_0 I^{-\alpha}\,, \forall I$,
 \item[b)] $d_I \le d_0 e^{-a I^{-\alpha}}\,, \forall I \Rightarrow \tilde\sigma_I \le C d_0 e^{-c I^{-\beta}}\,, \forall I $,
\item[c)] $d_I \le d_0 e^{-a I^{-\alpha}}\,, \forall I \Rightarrow \tilde\sigma_I \le C I \eta^I e^{-a I^{-\alpha}}\,, \forall I$,
\end{itemize}
see~\ref{app:weakgreedy} for a proof. Of course, the constants $C$ in the upper-bounds above can be overly pessimistic. But our numerical results show that one does not need to go to the asymptotics $I\to\infty$ in order to get good variance reduction results.

%\subsection{Implementation of the RB control-variate MC method}
\subsection{Implementation and analysis}

At any step $I\in\N$ of a weak greedy algorithm, whether we are looking for a new parameter value $\tilde\lambda_{I+1}$ or not ($TOL$ reached), the practical implementation of our RB control-variate MC method still requires a few discretization ingredients. In particular, one has to be able to actually compute approximations of $E(Z^\lambda)=E(Z^{\lambda}-\sum_{i=1}^I\alpha_i^\lambda Z^{\tilde\lambda_i})+\sum_{i=1}^I\alpha_i^\lambda E(Z^{\tilde\lambda_i})$ and $V(Z^\lambda-\hat Y^\lambda)$ whatever $\lambda\in\Lambda$. We are now going to suggest a generic practical approach to that question. On the contrary, practically useful choices of $\tilde\Lambda\subset\Lambda$ specifically depend on the parametric context and will be discussed only for specific applications in the next sections. 

First, we need to approximate the expectations $E(Z^{\tilde\lambda_i})$, $i=1\ldots I$. To obtain good useful approximations of $E(Z^{\tilde\lambda_i})$ in general, we suggest to use an expensive MC estimation $E_{M_{\rm large}}(Z^{\tilde\lambda_i})$, independent from subsequent realizations of the random variable $Z^{\tilde\lambda_i}$. Although this is a priori as computationally expensive as approximating well $E(Z^{\lambda})$ at any $\lambda$ by a naive MC estimation, we do it only once at each of the few $\tilde\lambda_i$. Once a realization of $E_{M_{\rm large}}(Z^{\tilde\lambda_i})$ has been computed at step $i$ of the weak greedy algorithm ($1\le i\le I$) with $M_{\rm large}\gg 1$ realizations, the deterministic result can be stored in memory. We denote the computable approximations of $Y^{\tilde\lambda_i}\equiv Z^{\tilde\lambda_i}-E(Z^{\tilde\lambda_i})$ by $\tilde Y^{\tilde\lambda_i}\equiv Z^{\tilde\lambda_i}-E_{M_{\rm large}}(Z^{\tilde\lambda_i})$, $i=1\ldots I$, so that the actual control variate $\hat Y^\lambda$ used in practice is in fact a linear combination in a linear space of dimension $I$ spanned by $\tilde Y^{\tilde\lambda_i}$, $i=1\ldots I$. Note that given $E_{M_{\rm large}}(Z^{\tilde\lambda_i})$ we expect realizations of $\tilde Y^{\tilde\lambda_i}$ to be computed concurrently with realizations of $Z^{\lambda}$, using the same pseudo-random numbers generated by a computer, with approximately the same cost. 

Second, for any $\lambda$, we need to replace the coefficients $\alpha_i^\lambda$, $i=1\ldots I$, with a tractable solution $\tilde \alpha^\lambda_i$ %to a computable approximation 
of the minimization problem~\eqref{eq:variance_minimization}, so that practical control variates in fact read
\beq
\label{eq:controlvariate}
\hat Y^\lambda = \sum_{i=1}^I \tilde \alpha_i^\lambda \tilde Y^{\tilde\lambda_i} \,.
\eeq
Now, the MC computations at any one specific $\lambda$ should be fast, in particular the approximate solution to the minimization problem~\eqref{eq:variance_minimization}, whether we decide to stop the weak greedy algorithm at step $I$ and use only $I$ parameter values or still want to explore the parameter values in $\tilde\Lambda$ to select a new parameter value $\tilde\lambda_{I+1}$. So we invoke only a {\it small} number $M_{\rm small}$ of realizations of $Z^{\lambda}$ and $\tilde Y^{\tilde\lambda_i}$ to compute the $\tilde \alpha^\lambda_i$.
Note that there are different strategies for the numerical solution to a least-squares problem like~\eqref{eq:variance_minimization} that use only $M_{\rm small}$ realizations of $Z^{\lambda}$ and $\tilde Y^{\tilde\lambda_i}$. One can try to compute directly the solution to an approximation of the linear system~\eqref{eq:RBnormal} where variances and covraiances are approximated with MC estimations similar to~\eqref{eq:varianceMCestimation} using $M_{\rm small}$ i.i.d. realizations $Z^{\lambda}_m$ and $\tilde Y^{\tilde\lambda_i}_m$, $m=1\ldots M_{\rm small}$, of $Z^{\lambda}$ and $\tilde Y^{\tilde\lambda_i}$. The MC estimations similar to~\eqref{eq:varianceMCestimation} are indeed interesting insofar as they remain positive, see~\cite{welford-1962,west-1979} for a study of consistency. But there may still be difficulties when $\mathbf{C}$ is too ill-conditioned. One can thus also apply a QR decomposition approach to a small set of $M_{\rm small}$ realizations, using the Modified Gram-Schmidt (MGS) algorithm for instance
\footnote{
  The matrix with entries $(\tilde Y^{\tilde\lambda_i}_m)_{m=1\ldots M_{\rm small},i=1\ldots I}$ uniquely decomposes as $QR$ where $Q$ is a $M\times I$ matrix such that $Q^TQ$ is 
  the $I\times I$ identity matrix and $R$ is an upper triangular $I\times I$ matrix.
  The $I$-dimensional vector $\tilde \alpha^\lambda$ with entries $\tilde \alpha^\lambda_{i}$, $i=1\ldots I$, useful in $\hat Y^{\lambda}$   is next defined as $R^{-1}Q^TZ$ (where $Z$ is the $M_{\rm small}$-dimensional vector with entries $Z^{\lambda}_m$). The $\tilde \alpha^\lambda_{i}$ are thus random variables depending on $M_{\rm small}$ realizations
  %but not on $M_{\rm large}$ realizations (through approximations $E_{M_{\rm large}}(Z^{\lambda_i})$ of $E(Z^{\lambda_i})$) !!	
  and an approximate minimum of~\eqref{eq:variance_minimization} is $Z^TZ-Z^TQQ^TZ$.
}.
In any case, realizations of $Z^{\lambda}$ and $\tilde Y^{\tilde\lambda_i}$, $i=1\ldots I$, should be correlated and thus computed with the same collection of random numbers in practice.
Besides, one can use a single collection of random numbers (generated by initializing the random number generator at one fixed seed) for all $\lambda$, including $\tilde\lambda_i$, $i=1\ldots I$, and realizations of $\tilde Y^{\tilde\lambda_i}$ can thus be computed only once (at step $i$ of the weak greedy algorithm) and next stored in memory.
%\footnote{
%  The MC estimations similar to~\eqref{eq:varianceMCestimation} are interesting insofar as they remain positive, see~\cite{welford-1962,west-1979} for a study of consistency.
%   The MC estimations should remain positive definite: compute iteratively 
%   $V_m\left(Z^{\lambda_i}\right)\approx V\left(Z^{\lambda_i}\right)$
%   and $E_m\left(Z^{\lambda_i}\right)\approx E\left(Z^{\lambda_i}\right)$ by
%   $$ 
%   E_{m+1}\left(Z^{\lambda_i}\right) = 
%   \frac1{m+1}\left( m E_m\left(Z^{\lambda_i}\right) + Z^{\lambda_i}_{m+1} \right) \,,
%   $$
%   $$ 
%   V_{m+1}\left(Z^{\lambda_i}\right) = 
%   \frac1{m+1}\left( m V_m\left(Z^{\lambda_i}\right) + |Z^{\lambda_i}_{m+1}-E_{m+1}(Z^{\lambda_i}|^2 \right) \,,
%   $$
%   thus in two steps for all $m=1\ldots M$.
%   But it is more prone to rounding error than computing sums and dividing them once only at the last step. 
%  see also welford-1962,west-1979,chan-golub-leveque-1979
%  But there may still be difficulties when $\mathbf{C}$ is too ill-conditioned.
%}

Finally we approximate $E(Z^{\lambda})$ by one draw of 
\beq
\label{eq:output0}
\frac1{M_{\rm test}} \sum_{m=1}^{M_{\rm test}} \left( Z^{\lambda}_m- \sum_{i=1}^I \tilde \alpha^\lambda_{i} \tilde Y^{\tilde\lambda_i}_m \right)
\eeq
whatever $\lambda$ and compute a MC estimation~\eqref{eq:varianceMCestimation} of $V(Z^{\lambda}-\hat Y^{\lambda})$ with the same $M_{\rm test}$ realizations of $Z^{\lambda}$ and $\tilde Y^{\tilde\lambda_i}$. Remember that at step $I\in\N$ of the weak greedy algorithm, \eqref{eq:output0} is useful either to inspect all parameter values $\lambda\in\tilde\Lambda$, check whether $TOL$ is reached and next select a new parameter value $\tilde\lambda_{I+1}$ when it is not. Or if $TOL$ has been reached before, then~\eqref{eq:output0} is used for definitive estimations at any $\lambda$ (still with companion variance estimations~\eqref{eq:varianceMCestimation} to concurrently certify that $TOL$ is maintained). 

Let us compare our strategy with the cost of direct MC estimations that have same confidence levels and thus use $M$ realizations such that $V_M(Z^{\lambda})/M = V_{M_{\rm test}}(Z^{\lambda}-\hat Y^{\lambda})/M_{\rm test}$. We denote by $\mathsf{C}$ the cost of one realization at one parameter value $\lambda$ compared with that of one multiplication. The evaluation of $M_{\rm small}+M_{\rm test}$ realizations at each $\lambda$, plus $I^2+IM_{\rm small}$ multiplications for the QR solution to the least-squares problem, $IM_{\rm test}$ multiplications for the MC estimation of the output expectation, and $M_{\rm test}^2$ for the variance, make our RB control-variate MC method interesting as soon as $ M\mathsf{C}+M^2 \ge (M_{\rm test}+M_{\rm small})(\mathsf{C}+I) + I^2 + M_{\rm test}^2$, at least for ``real-time'' purposes. Then, the price of identifying $I$ control variates with the weak greedy algorithm pays back if one needs to compute very fast $E(Z^\lambda)$ for any $\lambda$, provided the $I$ control variates still provide good variance reduction in the case where they were constructed by a weak greedy algorithm trained on $\tilde\lambda\subsetneq\Lambda$. In addition, the RB control-variate MC method is also interesting in the many-query cases where one has to compute $E(Z^\lambda)$ for a sufficiently large number $\sharp\lambda$ of parameter values $\lambda$. At each greedy step $i=1\ldots I-1$, in addition to variance estimations at each $\lambda\in\tilde\Lambda$, the selection of $\tilde\lambda_{i+1}$ requires a quicksort of the sample of estimated variances 
$\{V_{M_{\rm test}}\left( Z^{\lambda}-\hat Y^\lambda\right),\,\lambda\in\tilde\Lambda\}$, plus the computation of $E_{M_{\rm large}}(Z^{\tilde\lambda_{i+1}})$ and $M_{\rm small}$ realizations of $Z^{\tilde\lambda_{i+1}}$ to be stored. The cost of the greedy algorithm
$\mathsf{C}_{\rm greedy} = {\rm card}(\tilde\Lambda)\bigg(I(M_{\rm test}+M_{\rm small})\mathsf{C} + I(I+1)(2I+1)/6 + I(I+1)M_{\rm small}+M_{\rm test}/2+IM_{\rm test}^2+\ln{\rm card}(\tilde\Lambda)\bigg) + I(M_{\rm large}+M_{\rm small})\mathsf{C}$ % we do not need to evaluate expectations during greedy
is then compensated by variance reductions as soon as $\sharp\lambda \ge \mathsf{C}_{\rm greedy}/(M\mathsf{C}+M^2 - (M_{\rm test}+M_{\rm small})(\mathsf{C}+I) - I^2 - M_{\rm test}^2)$. 
% In the particular case where $\sharp\lambda={\rm card}(\tilde\Lambda)$ in particular, this requires
% \begin{multline*}
% M\mathsf{C}+M^2 - (M_{\rm test}+M_{\rm small})\mathsf{C} - I^2 - IM_{\rm small} - IM_{\rm test} - M_{\rm test}^2 \\
% \ge \mathsf{C}_{\rm greedy}/\sharp\lambda \approx I(M_{\rm test}+M_{\rm small})\mathsf{C} + I^3 + I^2M_{\rm small} + M_{\rm test} +IM_{\rm test}^2 + I(M_{\rm large}+M_{\rm small})\mathsf{C}/\sharp\lambda\,.
% \end{multline*}

%\subsection{Analysis of RB control-variate MC method}

We also have a posteriori error estimates. For values $\tilde \alpha^\lambda_{i}$ given by a fixed MC estimation with $M_{\rm small}$ realizations independent of $M_{\rm test}$ new realizations in~\eqref{eq:output0}, CLT~\eqref{eq:CLT} holds\footnote{
  If the $M_{\rm test}$ realizations in~\eqref{eq:output0} are the same as the $M_{\rm small}$ ones used to compute $\tilde \alpha^\lambda_{i}$, then~\eqref{eq:output0} is not a MC empirical estimator of the type $E_{M_{\rm test}}(Z^{\lambda}-\hat Y^{\lambda})$. It does not use independent realizations of the random variable $Z^{\lambda}-\hat Y^{\lambda}$, so the CLT does not hold 
  %\footnote{Chebyshev inequality still holds, but it may be difficult to estimate the variance.}
  and it is difficult to give a rigorous quantitative estimate of the statistical error.
}, $\forall a>0$, for 
\beq
%\begin{multline}
\label{eq:CLT2}
P\left( 
 % \left| E_M((Z^{\lambda}-Y^{\lambda})_m) - E(Z^{\lambda}-Y^{\lambda}) \right| \le a\sqrt{\frac{V(Z^{\lambda}-Y^{\lambda})}{M}} 
 \frac{\left| E_{M_{\rm test}}(Z^{\lambda}-\hat Y^{\lambda}) - E(Z^{\lambda})+ \sum_{i=1}^I \tilde\alpha^\lambda_{i}E(\tilde Y^i)) \right|}
      {\sqrt{V(Z^{\lambda}-\hat Y^{\lambda})/M_{\rm test}}} \le a 
\right)
%\\ 
%\underset{M\to\infty}{\longrightarrow} %\int_{-a}^a \frac{e^{-x^2/2}}{\sqrt{2\pi}} dx \equiv \mathop{\rm erf}\left(\frac{a}{\sqrt{2}}\right)
%\end{multline}
\eeq
where the MC empirical estimator $E_{M_{\rm test}}(Z^{\lambda}-\hat Y^{\lambda})$ is defined in~\eqref{eq:output0} and where one can replace $V(Z^{\lambda}-\hat Y^{\lambda})$ by a MC estimator $V_{M_{\rm test}}(Z^{\lambda}-\hat Y^{\lambda})$ by Slutsky theorem.
% the variance which also reads
% $$
% V(Z^{\lambda})+ \sum_{j=1}^I \tilde \alpha^\lambda_j \left( \sum_{i=1}^I \mathbf{C}_{i,j} % C(\tilde Y^j,\tilde Y^i)
% \tilde \alpha^\lambda_i -2 C(Z^{\lambda},Z^{\lambda_j}) \right) \,.
% $$
% \begin{equation}
% \label{eq:RBvarianceMsmall}
% % V(Z^\lambda-Y^\lambda) \approx 
% E_{M_{\rm small}}\left( \left| Z^{\lambda}-\tilde Y^{\lambda}-E_{M_{\rm small}}(Z^{\lambda}-\tilde Y^{\lambda})\right|^2 \right) \,.
% \end{equation}
Furthermore, although~\eqref{eq:CLT2} is not a full error analysis because it does not take into account the bias
$\sum_{i=1}^I \tilde\alpha^\lambda_{i}E(\tilde Y^i)=\sum_{i=1}^I \tilde\alpha^\lambda_{i}(E(Z^{\lambda_i})-E_{M_{\rm large}}(Z^{\lambda_i}))$, a function of $M_{\rm_{\rm large}}$ and $M_{\rm small}$ precomputed realizations, probabilities like~\eqref{eq:CLT2} are conditionally to $E_{M_{\rm large}}(Z^{\lambda_i})$ (and to $\tilde\alpha^\lambda_{i}$), so Bayes rule applies and it also holds, for all $a,a_i>0$:
\begin{multline}
\label{eq:CLT3}
P\bigg( 
 \left| E_{M_{\rm test}}(Z^{\lambda}-\hat Y^{\lambda}) - E(Z^{\lambda}) \right|
%      \\ \le a \sqrt{\frac{V_{M_{\rm test}}(Z^{\lambda}-\hat Y^{\lambda})}{M_{\rm test}}}
+ \sum_{i=1}^I |\tilde\alpha^\lambda_{i}| a_i \sqrt{\frac{V_{M_{\rm large}}(Z^{\lambda_i})}{M_{\rm large}}}
\bigg)
\\
\underset{M_{\rm large}\ge M_{\rm test}\to\infty}{\longrightarrow} % \int_{-a}^a \frac{e^{-x^2/2}}{\sqrt{2\pi}} dx \equiv 
\mathop{\rm erf}\left(\frac{a}{\sqrt{2}}\right) \prod_{i=1}^I \mathop{\rm erf}\left(\frac{a_i}{\sqrt{2}}\right)
\end{multline}
where variances have been replaced by estimators.

Last, to get a full convergence analysis that predicts an efficient variance reduction with the RB control-variate MC method in practice, at least when $\tilde\Lambda=\Lambda$, one should take into account all realizations in $\hat Y^\lambda$ which in fact reads $\hat Y^\lambda_{\rm M_{\rm small},M_{\rm large}} =\sum_{i=1}^I \tilde\alpha_{i,{\rm M_{\rm small}}}^\lambda \left( Z^{\tilde\lambda_i}-E_{M_{\rm large}}(Z^{\tilde\lambda_i}) \right)$. Yet, note first that in the weak greedy algorithm~\eqref{eq:weakgreedy}, only the $M_{\rm small}$ realizations introduce new statistical error. And second, one can again predict that, with a good probability, the greedy algorithm is robust to discretization (in the same sense as in~\cite{binev-cohen-dahmen-devore-petrova-wojtaszczyk-2010}) provided the realizations of the least-squares problems are not too close to rank-deficient and their numerical solution is close to the solution. We do not state this more rigorously but instead refer to~\cite{binev-cohen-dahmen-devore-petrova-wojtaszczyk-2010}, whose results can be adapted in the same way as in~\ref{app:weakgreedy} for the week greedy algorithm.

% \section{Application to Computation of Output after Uncertainty Propagation through Stochastic Coefficients in Partial Differential Equations}
\section{Application to Uncertainty Propagation}
\label{sec2}

In this section, we numerically demonstrate the efficiency of the RB control-variate MC method for uncertainty propagation in a representative UQ framework. As example, we consider a PDE parametrized by stochastic coefficients %that define a probability space 
and other non-stochastic coefficients which we term {\sl control parameters}. The goal is to compute fast many expectations of a scalar output of the random PDE solution for many values of the control parameters.

For the numerical simulations, as usual in UQ frameworks~\cite{ghanem-spanos-1991}, we use stochastic coefficients that are random fields with a Karhunen-Lo\`eve (KL) spectral decomposition. In practice, one can only use representations with a finite-rank $K$ of course and the MC method allows one to simulate the law of the random field just by generating realizations of the $K$ random numbers in the KL decomposition. We compute approximations to the realizations of the random PDE solution with a standard FE discretization method (the same one whatever the realizations). The expectation of the random scalar output can be computed for many values of a control parameter just by reiterating many times the MC procedure. But this is costly. So first, we do not actually compute the PDE solution with the FE method at each realization of the stochastic coefficients and for each control parameter value. In fact, we replace it with a cheap reliable surrogate built with the standard RB method like in~\cite{boyaval-lebris-maday-nguyen-patera-2009}. Yet, even with a cheap surrogate model instead of the full FE, the MC method is still costly, because computing the expectation of the parametrized output under the input probability law still requires a very large number of realizations. This quickly becomes untractable as soon as one has to do it many times for many values of the control parameter. Then, we use the RB control-variate MC method to decrease the number of MC realizations needed at most of the control parameter values.

Compared with~\cite{boyaval-lebris-maday-nguyen-patera-2009}, one uses here an improved MC method to show how various RB ideas can be combined to efficiently tackle some uncertainty propagation problems for partial differential equations with stochastic coefficients. We thus especially discuss the practice of the RB control-variate MC method here. Yet, although the example is exactly the same as in~\cite{boyaval-lebris-maday-nguyen-patera-2009}, the specific use of the standard RB method in the frame of UQ, and in combination with the RB control-variate MC method, requires special care. That is why at the end of this section we briefly expose our implementation of the standard RB method, without which we cannot use the MC method at a reasonable price, and next discuss our specific use of it.

\subsection{An elliptic PDE with stochastic coefficients}
\label{sec:parametrizedPDE}

\mycomment{RB for unregular geometries, unsymmetric/random $\bm{A}$, mixed problems (Dirichlet-Robin/Neumann) ?}

Consider first the solution $u$ to a scalar Robin Boundary Value Problem (BVP) in a regular domain $\D\subset\R^2$: 
$u$ satisfies Laplace equation in $\D$ % (simply-connected bounded Lipschitz polyhedral)
\begin{equation} 
\label{eq:strong}
 - {\rm div} \left( \kappa \nabla u \right) = 0 \,, % \bm{A}
\end{equation}
and conditions of third type on the boundary $\partial \D$ %smooth one-dimensional manifold (with finite one-dimensional measure)
\begin{equation}
\label{eq:BC}
% \bm{n}\cdot \bm{A} \nabla u + b \: u = g \,,
\kappa (\bm{n}\cdot\nabla) u + b \: u = g \,,
\end{equation}
where
$\rm div$ and $\nabla$ %respectively 
denote the usual divergence and gradient operators
in $\D$ equipped with a cartesian frame,
% $\bm{A}$ is a 2-tensor field in $\D$
% with %positive-definite %symmetric
% isotropic diagonal matrix values
% $$
% \bm{A} = \begin{pmatrix} k & 0 \\ 0 & k \end{pmatrix}
% $$
% while $\bm{n}\cdot \bm{A}$ denotes its contraction with
$\bm{n}$ the outward unit normal on $\partial \D$, and $k$, $b$ and $g$ are scalar parameter functions. For the simulations we will choose in particular
\begin{align}
\kappa & = k_1 \mathrm{1}_{\D_1} + k_2 \mathrm{1}_{\D_2} && (k_1,k_2)\in\R_{>0}^2
\\
b & = b \mathrm{1}_{\Gamma_{\rm B}} && b\in L^\infty(\partial\D,\R_{>0})
\\
g & = g \mathrm{1}_{\Gamma_{\rm R}} && g\in L^2(\partial\D)
\end{align}
%where $\mathrm{1}_{\D_i}$ denotes the characteristic function of Lipschitz subdomains $\D_i$ ($i=1,2$)
where $\mathrm{1}_{\D_i}$ is the characteristic function of $\D_i$,
$$
\D_1 \cap \D_2 = \emptyset
\quad
\D_1 \cup \D_2 \subset \D \subset \overline{\D_1 \cup \D_2}
\,,
$$
and the boundary decomposes into % non-overlapping 
subsets
$$
\Gamma_{\rm B} \cap \Gamma_{\rm R} = \emptyset
\quad
\Gamma_{\rm B} \cup \Gamma_{\rm R} \subset \partial\D
\quad
\Gamma_N =  \partial\D \backslash \overline{\Gamma_{\rm B} \cup \Gamma_{\rm R}} \,.
$$
There exists a unique weak solution $u \in H^1(\D)$ to~(\ref{eq:strong}--\ref{eq:BC}), which can also be defined as the unique solution to the following variational problem~\cite{quarteroni-valli-1997}: Find $u \in H^1(\D)\,/\,\forall v \in H^1(\D)$
\begin{equation}
\label{eq:weak}
% k_1 \int_{\D_1} \nabla u \cdot \nabla v + k_2 \int_{\D_2} \nabla u \cdot \nabla v 
\int_{\D_1} k_1 \nabla u \cdot \nabla v + \int_{\D_2} k_2 \nabla u \cdot \nabla v 
+ \int_{\Gamma_{\rm B}} b \: u \: v = \int_{\Gamma_{\rm R}} g \: v \,. %\,,\ \forall v \in H^1(\D) \,.
\end{equation} 
%The trace of $u,v \in H^1(\D)$ is well-defined on $\partial\D$.
%that problem is symmetric % ($\bm{A} =\bm{A} ^t$)
As output, we consider the compliance
\footnote{
  For symmetric problems like~\eqref{eq:weak}, the output~\eqref{eq:output} is a particularly simple choice because it allows a simple accurate a posteriori error estimation without invoking a dual problem, but neither this choice nor the symmetric character of the problem are limitations to our approach, see e.g.~\cite{nguyen-veroy-patera-2005}.
}
\begin{equation}
\label{eq:output}
s \equiv l(u) = \int_{\Gamma_{\rm R}} g \: u \,.
\end{equation}
The accurate numerical discretization of~\eqref{eq:weak} is standard, for instance using the Finite-Element (FE) method. In what follows, we shall invoke
continuous piecewise linear approximations defined in conforming, regular FE spaces of $H^1(\D)$.
% (endowed with the usual approximation properties  necessary for the convergence $u_{\cal N}(\mu)\to u(\mu)$
% to be uniform in $\mu$ as the number $\cal N$ of degrees of freedom increases).

\begin{figure}
\centering
%  \begin{tabular}{cc}
è    \makebox[.5\textwidth][c]{\includegraphics[scale=.5]{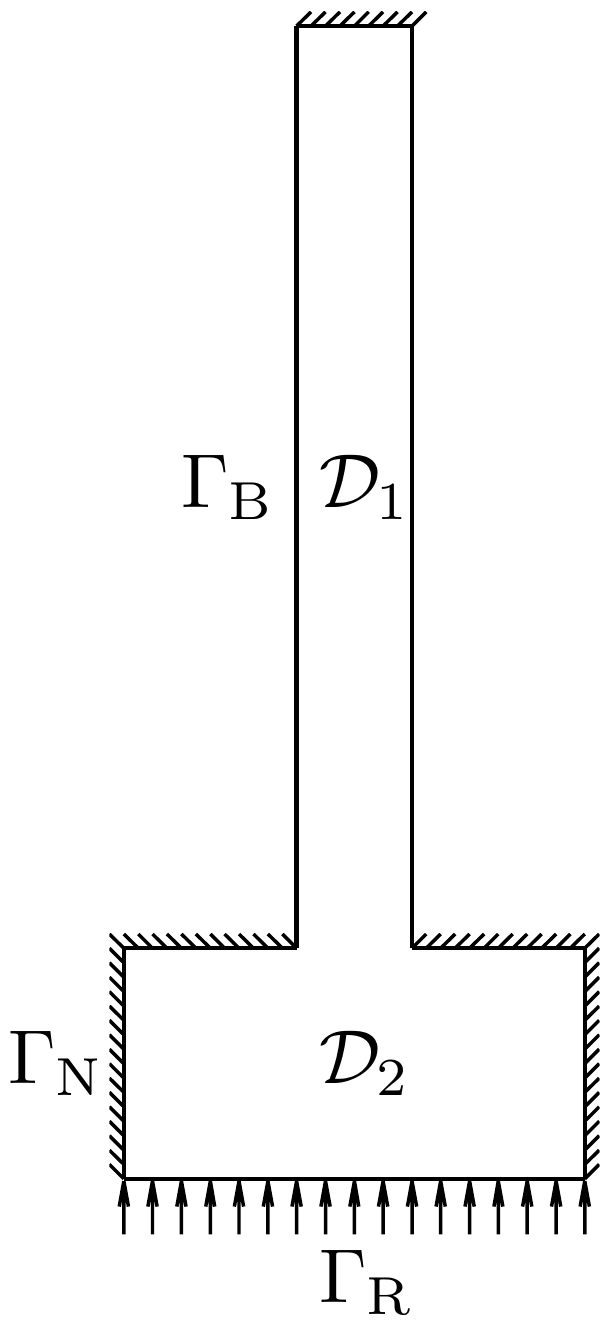}}
%    \makebox[.3\textwidth][c]{\includegraphics[scale=.45]{fig/Th_fin.pdf}}
%  \end{tabular}
\caption{\label{fig:geometry}
Thermal fin geometry $\D$ used in the numerical simulations
($\Gamma_{\rm N} = {\overline{\partial\D \backslash (\Gamma_{\rm R}\cup\Gamma_{\rm B})}}$ denotes the boundary subject to homogeneous pure Neumann conditions).
}
\end{figure}

%\subsection{Problem parametrization and discretization}\label{sec:1-prop}

\mycomment{ Vector stochastic coefficient fields: how to optimally represent them in a separate fashion, dealing with possible correlations between components ? 
For well-posedness, treat also log-normal coefficients with Empirical Interpolation. Can we hope for fast assembling (very few magic points) even with large $K$ ?}

To fix ideas, one could think of $u$ as the temperature in a fin subject to a constant radiative flux $g$ on $\Gamma_{\rm R}$ (in contact with a heat source)
and to convective thermal exchanges on $\Gamma_{\rm B}$, see Fig.~\ref{fig:geometry}. The function $b$ determines the value of the Biot number along $\Gamma_{\rm B}$,
that is the intensity of local heat transfers on the part of the boundary in contact with a fluid. %convective 
Now, the Biot number is typically uncertain. 
% \footnote{
%   Let us note that it makes sense to let $b$ be an uncertain (random) coefficient because quite often, the Biot number only roughly models the convective heat transfers
%   in the physical problems modelled of the Robin problem~(\ref{eq:strong}--\ref{eq:BC}).
%   And it is natural to evaluate the average of fluctuations induced by the stochastic variations of $b$ on quantities of interest like $s$ as a function of another parameter like $k_2$ for engineering purposes.
% }.
We shall model it as a random field in a probability space $(\Omega,\mathcal{F},P)$ and one is typically interested in quantifying the uncertainty propagated in $s$. More precisely, we define an essentially bounded function $b\in L^\infty_P(\Omega,L^\infty(\Gamma_{\rm B},\R_{>0}))$, positive almost everywhere (a.e.) on $\Gamma_{\rm B}$ with probability $1$ ($P$-almost surely or a.s. in abbrev.). Then, with a slight abuse of notation, $u$ now denotes the solution in $L^\infty_P(\Omega,H^1(\D))$ to~(\ref{eq:strong}--\ref{eq:BC}) when $b\in L^\infty_P(\Omega,L^\infty(\Gamma_{\rm B},\R_{>0}))$. And, for many values of the control parameter $k_2$, we consider computing expectations $E(s) = \int_{\Omega} s dP $ under $P$ of the output $s$. To that aim, we will propagate the uncertainty (random ``noise'') from $b$ to the scalar (random) output variable $s\in L^\infty_P(\Omega)$ by simulating the law of $u$.
% (It is easy to give a meaning to this since the deterministic $u\in H^1(\D)$ is continuous with respect to $b\in L^\infty(\Gamma_{\rm B},\R_{>0})$ and
% any $b$ in the Banach space $L^\infty_P(\Omega, L^\infty(\Gamma_{\rm B},\R_{>0}))$ can be approximated as close as wanted using finitely many random numbers.)
% Although a $L^2_P(\Omega)$ discretization framework (where approximations of $s$ converge to $s$ only in the $L^2_P(\Omega)$ sense)
% suffices for $\E{s}$, uncertainty quantification also often consists in providing one with other information, 
% for instance the probabilities of some events like $P(\{s>s_0\})$ (think of risk assessment), and thus
We recall that here we focus on % a $L^\infty_P(\Omega)$ discretization framework of
MC discretizations of the noise. Compared with other discretizations~\cite{matthies-keese-2004,babuska-nobile-tempone-2007}, it is very easy to implement, in particular because it is not intrusive with respect to existing discretizations of the deterministic problem, although the convergence is a priori much slower when the stochastic variations are smooth~\cite{cohen-devore-schwab-2010,cohen-devore-schwab-2011}.

%But first, we shall precise the law of the random input $b$ that we consider for numerical simulations.
% In spite of the MC $L^\infty_P(\Omega)$ framework, we nevertheless exploit the Hilbert structure of $L^2_P(\Omega)$ topology 
% to practically represent the random fluctuations induced by $b$ in an explicit, constructive way that separates (``tensorizes'') the dependences on $\Omega$ and $\D$.
We consider a bounded second-order stochastic process $b$ as random input. We denote its constant mean value on $\Gamma_{\rm B}$ by $\E{b|_{\Gamma_{\rm B}}}\equiv\bar E>0$, the spatial correlation length between random fluctuations by $\delta>0$ and the relative intensity of the (bounded) random fluctuations around the mean by $\Upsilon>0$. More precisely, like in~\cite{boyaval-lebris-maday-nguyen-patera-2009}, the random field $b$ shall be the limit in $L^\infty_P(\Omega,L^\infty(\Gamma_{\rm B}))$ of random fields $b_K$ whose covariance defines a kernel operator of finite rank $K$ in $L^2(\Gamma_{\rm B})$ and which can be easily simulated with a MC method using $K$ {\it independent} random scalar parameters in {\it a priori} given ranges. %This is sufficient for many applications of the MC method.
% Finite-rank kernels are indeed dense in the Banach spaces of Hilbert-Schmidt and trace-class operators % equipped with their respective norm
% Provided one can generate many pseudo-random realizations of $b$ and solve~(\ref{eq:strong}--\ref{eq:BC}) numerically for each value of $k_2$ and $b$, one can easily simulate the law of the scalar output $s$ of the random Robin problem~(\ref{eq:strong}--\ref{eq:BC}) at various values of $k_2$ and $\bar E$. And given a MC numerical simulation of the random output $s(\lambda)$ at a fixed $\lambda=(k_2,\bar E)$, one can easily compute MC estimations of the output expectation. (We will show how to use the RB control-variate MC method to improve the convergence rate of the MC estimations with respect to the number of realizations of $b$ in the limit of many queries for $\lambda$ only in the next section.)
For the numerical simulations, we choose a gaussian kernel $c(x,y)=\exp\left(-|x-y|^2/\delta^2\right)/|\Gamma_{\rm B}|$ to define a covariance operator in $L^2(\Gamma_{\rm B})$, where $x,y$ are in the same connected component of the boundary $\Gamma_{\rm B}$ in Fig.~\ref{fig:geometry}. The covariance operator has eigencouples $({\lambda_k},\Phi_k)$, $k\in\mathbb{N}_{>0}$ and a normalized trace $\sum^{\infty}_{k = 1} \lambda_k=1$. Then, invoking uniformly distributed independent variates $Z_k \sim \mathcal{U}(-\sqrt{3},\sqrt{3})$ for each spectral mode $k\in\mathbb{N}_{>0}$ (the same on each size of the boundary $\Gamma_{\rm B}$ in Fig.~\ref{fig:geometry} for the sake of simplicity), we define $b$ by its KL representation
\begin{equation} 
\label{KL}
b := \bar E \left( 1 + \Upsilon \sum_{k\in\mathbb{N}_{>0}} \sqrt{\lambda_{k}} \: \Phi_{k} \: Z_k \right) \,,
\end{equation}
%its KL representation, 
which makes sense since truncated representations
\begin{equation} 
\label{KL_K}
b_K = \bar E \left( 1 + \Upsilon \sum^{K}_{k = 1} \sqrt{\lambda_k} \: \Phi_{k} \: Z_k \right)
\end{equation}
converge in $L^2_P(\Omega,L^2(\Gamma_{\rm B}))$, and in $L^\infty_P(\Omega,L^\infty(\Gamma_{\rm B}))$ here (see~\cite{frauenfelder-schwab-todor-2005,schwab-todor-2006} e.g.).
We also require $b_K \ge b_{min} := \bar E/2 >0$ a.s. $\forall K\in\N$ and fix $\Upsilon:=.1$ for $\delta\in(.05,.5)$.
% \footnote{
%  In practice, we precompute good approximations of $({\lambda_k},\Phi_k)$ up to some order $K=K_{\rm max}$ by solving large eigenvalue problems, which gives an upper-bound for $\Upsilon$.
% }
% \begin{equation}
% \label{req:unif_bound}
% % b_{min} = \bar E \left( 1 - \Upsilon \sum_{k\in\mathbb{N}_{>0}} \sqrt{3\lambda_{k}} \|\Phi_{k}\|_{L^\infty(\Gamma_{\rm B})} \right) \,.
% 1/2 = \Upsilon \sum_{k\in\mathbb{N}_{>0}} \sqrt{3\lambda_{k}} \|\Phi_{k}\|_{L^\infty(\Gamma_{\rm B})} \,.
% \end{equation}
% (See the ten first eigenvectors in Fig.~\ref{fig:eigenvectors}, and observe the fast decay of the spectrum in Fig.~\ref{fig:variancebound}.)
% \begin{figure}
% \centering
% %  \begin{tabular}{cc}
% %    \makebox[.5\textwidth][c]{\includegraphics[scale=.4,trim= 0 60 0 100]{fig/eigenvectors1.pdf}}
%     \makebox[.35\textwidth][c]{\includegraphics[scale=.42,trim= 0 60 0 100]{fig/eigenvectors2.pdf}}
% %  \end{tabular}
% \caption{\label{fig:eigenvectors} %5 and 
% 10 first eigenvectors of the autocovariance operator. % computed using FEM.
% }
% \end{figure}
Then %, since $b_{K}>0$, 
one can define well-posed BVPs: %when $b$ is approximated by $b_{K}$: \\
Find $u_{K}\in L^\infty_P(\Omega,H^1(\D))$ such that $P$-a.s.
\begin{multline}
\label{eq:strongK}
\kappa (\bm{n}\cdot\nabla) u_K + b_K \: u_K = g \text{ on } \partial \D \,,
\\
- {\rm div} \left( \kappa \nabla u_K \right) = 0 \text{ on } \D \,.
\end{multline}
Approximations $u_K$ converge to $u$ in $L^\infty_P(\Omega,H^1(\D))$ as $K\to\infty$
and for a.e. $\omega\in\Omega$, realizations $u_K(\omega) \in H^1(\D)$ of the solution to~(\ref{eq:strongK})
solve a variational formulation with test function $v \in H^1(\D)$
\begin{multline}
\label{eq:weak_KL}
k_1 \int_{\D_1} \nabla u_K(\omega) \cdot \nabla v + k_2 \int_{\D_2} \nabla u_K(\omega) \cdot \nabla v 
\\
+ \bar E \sum_{k=1}^K (\Upsilon \sqrt{\lambda_k} Z_k(\omega)) \int_{\Gamma_{\rm B}} \Phi_k \: u_K(\omega) \: v 
\\
+ \bar E \int_{\Gamma_{\rm B}} e \: u_K(\omega) \: v
= \int_{\Gamma_{\rm R}} g \: v % \,,\ \forall v \in H^1(\D) 
\end{multline}
parametrized by $k_1,k_2,\bar E$ and realizations $\Upsilon \sqrt{\lambda_k} Z_k(\omega)$, $k=1\ldots K$. So finally, in practice, we shall compute approximations $u_{\mathcal{N},K}$ to $u$ with the FE method defined previously, which converge in $L^\infty_P(\Omega,H^1(\D))$ as $\mathcal{N}\to\infty$, $K\to\infty$. % H^1 approximation property does not depend on K
We require $\mathcal{N}$ and $K$ large enough such that, for a.e. values of the control parameter $k_1,k_2,\bar E$,
% the FE approximations $u_{\mathcal{N},K}$ to be uniformly good in $H^1(\D)$
a tolerance level $|s-s_{\mathcal{N},K}|\le TOL |s|$ is reached a.s. for the approximate output $s_{\mathcal{N},K}:=\int_{\Gamma_{\rm B}}u_{\mathcal{N},K}$.
Next one can easily simulate the law of $s_{\mathcal{N},K}$ with the MC method by mapping realizations $Z_k(\omega)$ of $Z_k$, $k=1\ldots K$, to realizations $s_{\mathcal{N},K}(\omega)$ of $s_{\mathcal{N},K}$.

But a direct MC-FE method is very expensive computationally, since one has to compute many FE approximations $u_{\mathcal{N},K}(\omega)$,
for many realizations of $Z_k(\omega)$, and next for many values of $\lambda=(k_2,\bar E)$.
That is why, like in~\cite{boyaval-lebris-maday-nguyen-patera-2009}, we replace the FE approximations $u_{\mathcal{N},K}$
by cheaper surrogates $u_{\mathcal{N},K,N}$ constructed with the standard RB method (wee will come back to this point in Section~\ref{sec:certifiedRBmethod} and the non-expert reader can find a brief exposition of this technique there).
% More generally, the BVP for $u$, parametrized by $k_1,k_2,b$ and $g$, has to be solved many times for many values of $(k_1,k_2,b,g)$ in a number of applications.
% It is often neither possible nor efficient to repeatedly deploy the FE method for each value of the (infinite-dimensional) parameter $b$. 
% We will use the RB method to compute fast the many FE approximations of the solutions to~\eqref{eq:weak} at many values of the parameter $(k_1,k_2,b,g)$. 
At present, the RB method is one of the only few existing alternatives that can tackle some ``high-dimensional'' problems in more than 2 or 3 dimensions, and it proved useful\footnote{
 Note that it does not a priori require specific discretizations of the parameter space and speed-up can be increased for goal-oriented purposes like computing $s$ accurately instead of $u$.
} when the parameter is not too high-dimensional, for instance when each of the parameter components $k_1,k_2,b$ and $g$ are scalars. But when the parameter becomes high-dimensional, difficulties arise again at some point. To some extent, the RB method still showed efficient for MC simulations of a random field $b$ with a moderately large number $K$ of modes~\cite{boyaval-lebris-maday-nguyen-patera-2009}. But there are situations where the approach is still computationally too expensive, in particular when one wants to explore variations of the MC simulation with respect to the control parameter $k_2$ and $\bar E$ for instance. That is why, in this work, we would next like to show how to further improve the computational cost in cases where one is interested in computing many values of the expectation $E(s)$ of $s$ under the law of the random field $b$ for many values of the control parameters $k_2$ and $\bar E$. To this aim, we build on the standard RB ideas and use a MC method with a reduced basis of control variates. In particular, at a given value of $\lambda$, one still needs a large number of RB approximations $u_{\mathcal{N},K,N}$ to compute accurate MC estimations of $\E{s_{\mathcal{N},K,N}}$, which is too costly when it has to be done for many values of $\lambda$. We thus now try to reduce the computational cost by invoking still another computational reduction technique, the RB control-variate MC method. We will show that the number of realizations required by the MC method to reach a given statistical accuracy for $\E{s_{\mathcal{N},K,N}}$ at one parameter value $\lambda$ can be efficiently decreased with our RB control-variate technique,
in the limit where one has to compute sufficiently many expectations for sufficiently many values of $\lambda$.

\subsection{RB control-variate MC estimations}

We require a relative precision $TOL=10^{-3}$ %of three digits 
for the numerical approximations $E(s_{\mathcal{N},K,N})(\lambda)$ of the output expectations $E(s)(\lambda)$ at any $\lambda\in\Lambda:=(.1,10)\times(.1,1)$. 
%Since $s$ is a compliant PDE output and the coercivity constant is bounded below by $.05$, 
It is thus sufficient for the PDE solution $u$ to be a.s. approximated with a relative precision $tol=10^{-2}$ here. To this aim, we first define $\mathbb{P}_1$ FE approximations on a 2D simplicial mesh using FreeFem++ by F.~Hecht and his collaborators. (Refining a coarse mesh, $\mathcal{N} = 3333$ nodes -- $6272$ triangles -- seem sufficient.) % The mesh is endowed with the symmetry of the geometry $\D$: a vertical symmetry line in the middle of the fin, like the random field $b$.

\begin{figure}
%\centering
  \begin{tabular}{r}
%  \begin{tabular}{c}
%    \makebox[.4\textwidth][c]{\includegraphics[scale=.31,trim= 50 95 10 95,clip]{fig/meanKL3N5.pdf}}
    \makebox[.4\textwidth][c]{\includegraphics[scale=.5,trim= 10 20 0 20,clip]{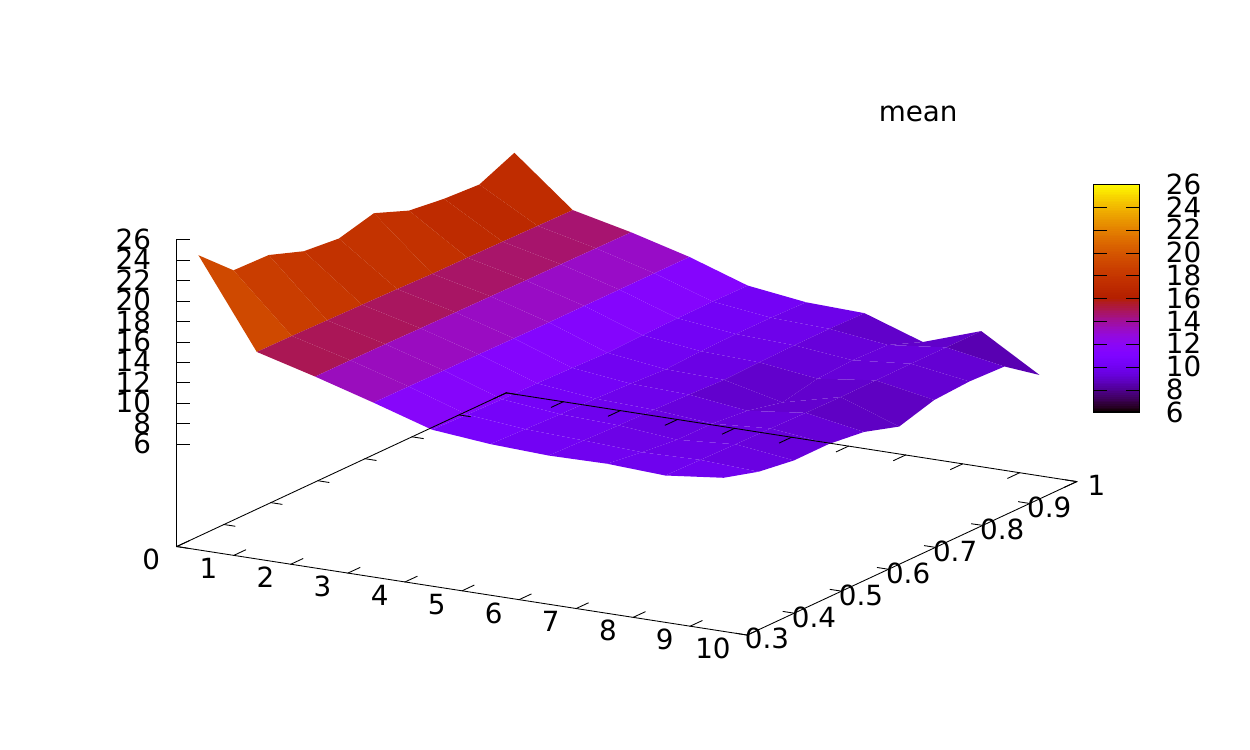}}
%  \end{tabular}
\\
% &  \begin{tabular}{c}
%    \makebox[.46\textwidth][c]{\includegraphics[scale=.31,angle=-90,trim= 10 10 40 10,clip]{fig/varKL3N5.pdf}}
    \makebox[.4\textwidth][c]{\includegraphics[scale=.5,trim= 10 20 0 20,clip]{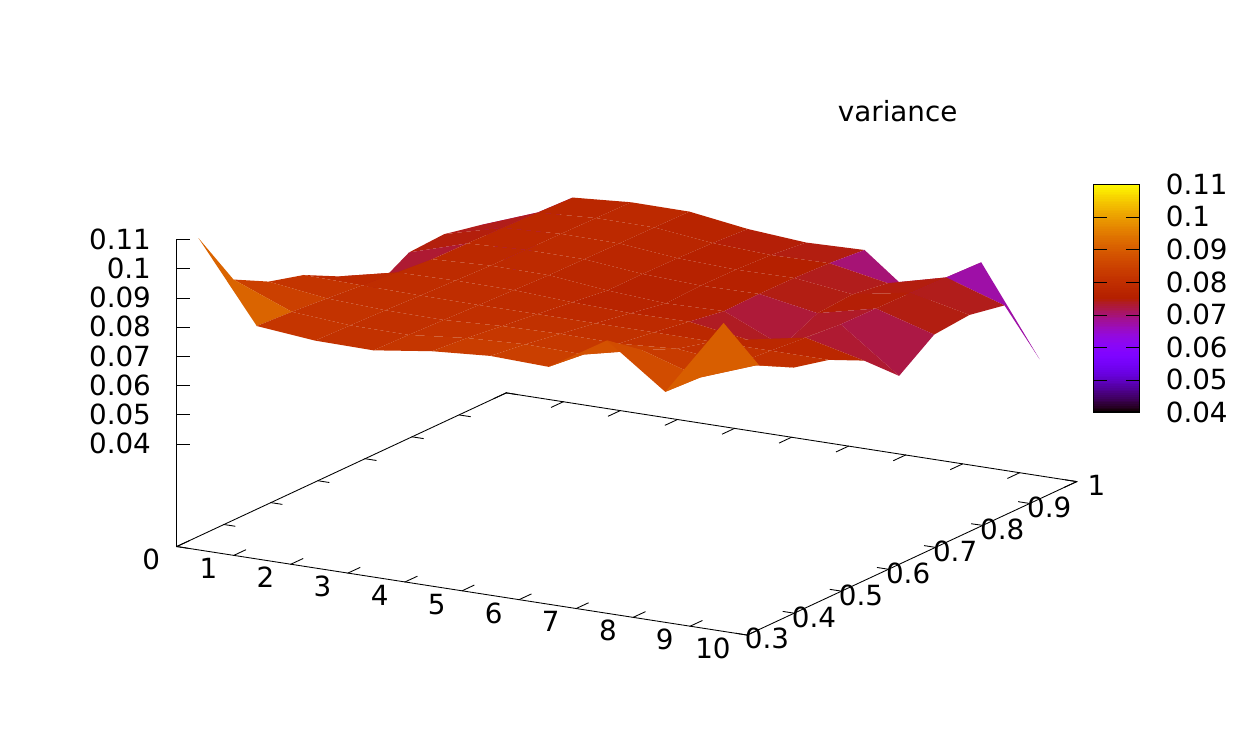}}
%   \end{tabular}
  \end{tabular}
\caption{\label{fig:MCmaps} Maps of $\E{s(k_2,\bar E)}$ and $\Var{s(k_2,\bar E)}$ (numerical approximations) for $(k_2,\bar E)\in(.1,10)\times(.1,1)$.
}
\end{figure}

We first fix $\delta = 0.5$. By a standard ``variational crime'' analysis (see e.g.~\cite{boyaval-lebris-maday-nguyen-patera-2009}), it is enough to use a truncated KL expansion $b_K$ up to $K=10$ (such that $\sum_{k>K}\sqrt{\lambda_k}\le tol\sum_{k\in\mathbb{N}_{>0}}\sqrt{\lambda_k}$). We build ``offline'' (that is, before the MC simulation) a reduced basis for a small linear space\footnote{
  A weak greedy algorithm needed $N=12$ steps for the a posteriori error estimates to satisfy $\|u_{\mathcal{N},K}-u_{\mathcal{N},K,N}\|_{H^1}\le\Delta_N\le tol=10^{-2}$ at all training parameter values in a cartesian grid using 2 trial values for each parameter component $\Upsilon \sqrt{\lambda_k} Z_k$, $k=1\ldots K$, and 10 trial values for each parameter component $\lambda_1=k_2,\lambda_2=\bar E$, see Section~\ref{sec:certifiedRBmethod} for some details about the RB method.
} of surrogates $u_{\mathcal{N},K,N}$ with dimension $N=12$. Then we simulate the law of $Z^\lambda\equiv s_{\mathcal{N},K,N}(\lambda)$ at any $\lambda\in\Lambda$ with the MC method. This allows us to retrieve a MC estimation of $E(Z^\lambda)$ and to construct a response surface by piecewise linear interpolation of the MC estimations $E_M(Z^\lambda)$ at the $10\times 10$ trial values of $\lambda$ like in Fig.~\ref{fig:MCmaps}. 

To compute the response surface of Fig.~\ref{fig:MCmaps} with a good confidence level in the MC estimations, say with $\sqrt{V_M(Z^\lambda)/M}\le TOL\times E_M(Z^\lambda)$ at each of the $\lambda$ used for interpolation, one needs more than $M=10^4$ realizations\footnote{
  Besides, one can check that $M=10^4$ is sufficiently-many realizations here for the speed of convergence to be quite well evaluated by CLT,
  since the MC estimator $E_M(Z^\lambda)$ is already close to gaussian according to the Kolmogorov-Smirnov goodness-of-fit test.
  In particular, most often, one single draw of the MC estimator already allows one to compute a good approximation to a confidence interval.
} at each $\lambda$. % $2.58\sqrt{V_M(Z^\lambda)/M}$ as approximate error bound with probability  99\%
Now, the law of $Z^\lambda$ with respect to $\lambda$ is smooth. Let us use the RB control-variate MC method with $\tilde\Lambda$ defined by the $10\times 10$ trial values of $\lambda$ needed to build the response surfaces in Fig.~\ref{fig:MCmaps}. Then, we can achieve MC estimation with confidence intervals of accuracy $TOL$ using a large number of realizations $M_{\rm large}=10^4$ at only $I=3$ values of $\lambda$, and much less realizations at the other values of $\lambda$ to achieve satisfying confidence intervals. More precisely, $M_{\rm test}=10$ is enough, since one can reduce the variance to $10^{-4}$ with $I=3$, and one already has a good MC estimation of the magnitude of the variance then. One also needs to compute coefficients $\alpha_i^\lambda$ %in the linear combination 
at each $\lambda$ with $M_{\rm small}$ additional realizations, but in practice one can use the same $M_{\rm test}=M_{\rm small}$ realizations\footnote{
 Although introducing some dependency between the $M_{\rm test}$ and $M_{\rm small}$ realizations theoretically influences the possibilty to obtain a CLT for instance, it hardly changes siginificantly the numerical result for one MC estimation in practice.
}.
So the marginal gain for each $\lambda$ in the finite sample of queries $\tilde\Lambda$ is to divide the computational cost by $10^4/(10+I\times10^4/Card(\tilde\Lambda))\ge30$
(the cost $\sf C$ of one realization is dominant here, even with RB surrogates).

\begin{figure}
%\centering
  \begin{tabular}{r}
%  \begin{tabular}{c}
%    \makebox[.4\textwidth][c]{\includegraphics[scale=.31,trim= 50 95 10 95,clip]{fig/meanKL3N5.pdf}}
    \makebox[.5\textwidth][c]{\includegraphics[scale=.5,trim= 10 20 0 20,clip]{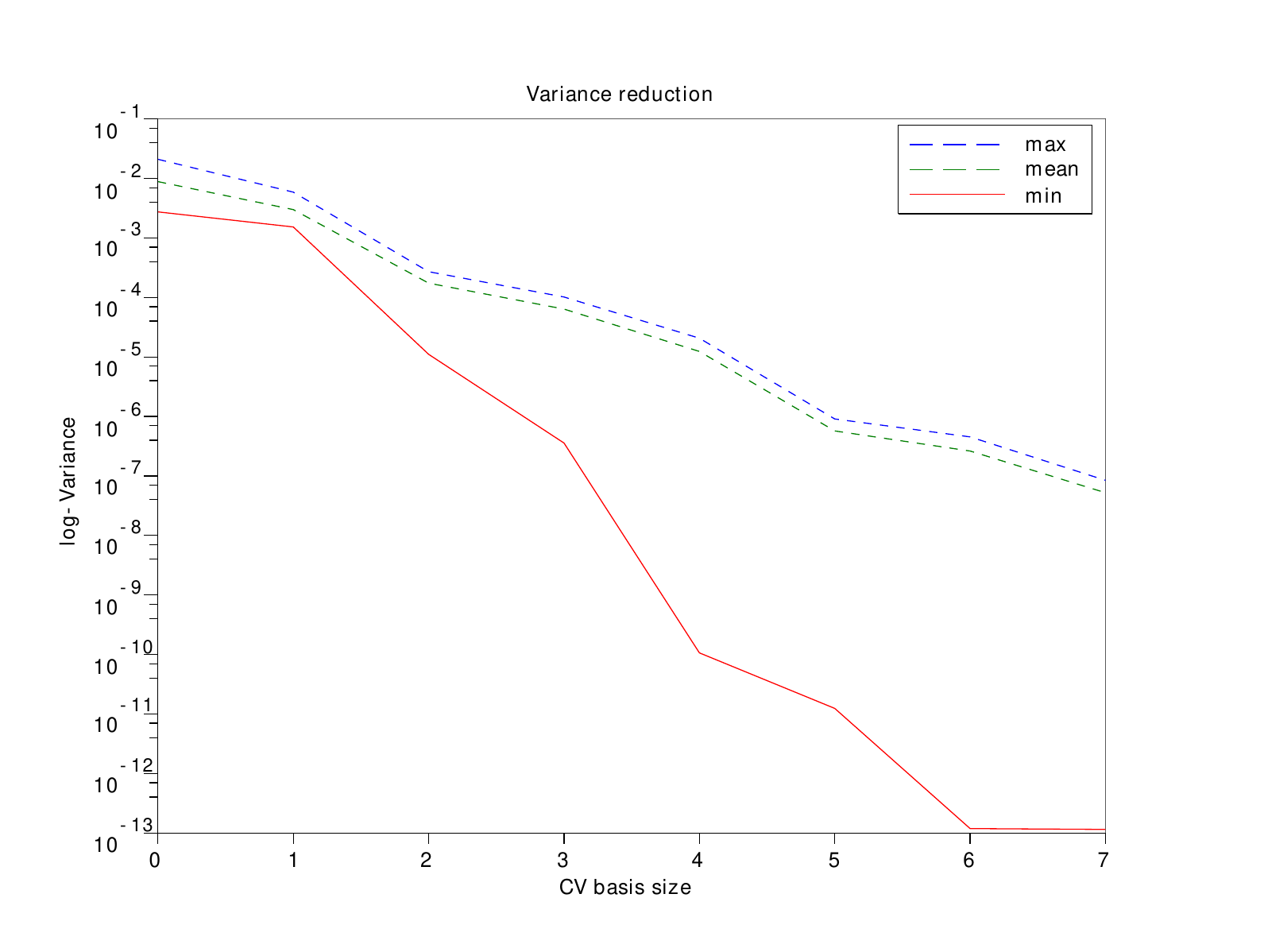}}
%  \end{tabular}
\\
% &  \begin{tabular}{c}
%    \makebox[.46\textwidth][c]{\includegraphics[scale=.31,angle=-90,trim= 10 10 40 10,clip]{fig/varKL3N5.pdf}}
%    \makebox[.4\textwidth][c]{\includegraphics[scale=.5,trim= 10 20 0 20,clip]{fig/vard5em1K10.pdf}}
%   \end{tabular}
%   \end{tabular}
% \caption{\label{fig:MCchecks} 
% Statistics among 10 new $\lambda$ of the empirical means of the estimated variance $V_{M_{\rm test}}(Z^\lambda-\sum_{i=1}^\iota\alpha_i^\lambda\tilde Y^i)$
% using 100 draws using $\iota=0\ldots I$ control variates (log-scale).
% }
%\end{figure}
%\begin{figure}
%\centering
%  \begin{tabular}{r}
%  \begin{tabular}{c}
%    \makebox[.4\textwidth][c]{\includegraphics[scale=.31,trim= 50 95 10 95,clip]{fig/meanKL3N5.pdf}}
    \makebox[.5\textwidth][c]{\includegraphics[scale=.5,trim= 10 20 0 20,clip]{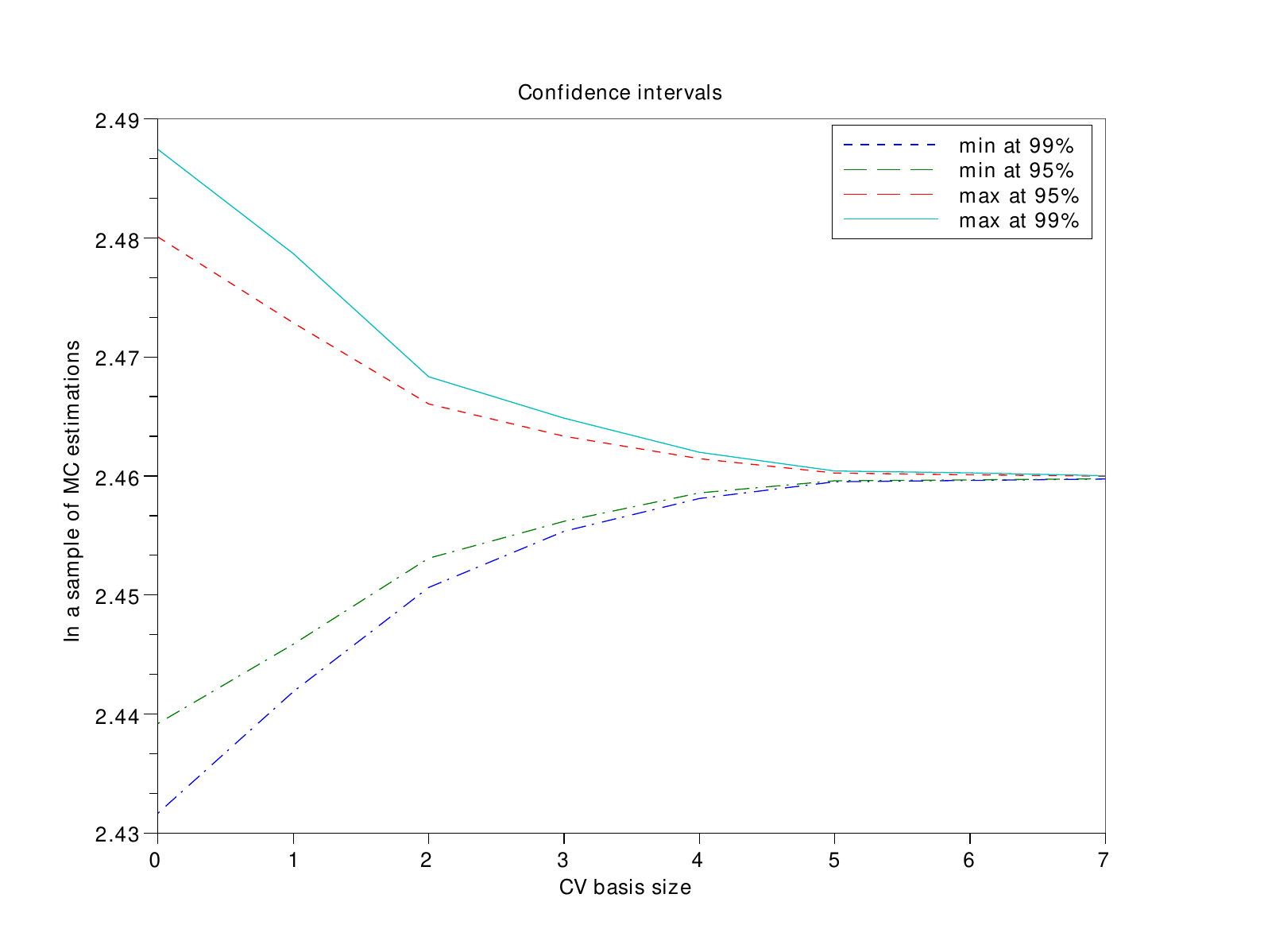}}
%  \end{tabular}
\\
% &  \begin{tabular}{c}
%    \makebox[.46\textwidth][c]{\includegraphics[scale=.31,angle=-90,trim= 10 10 40 10,clip]{fig/varKL3N5.pdf}}
    \makebox[.5\textwidth][c]{\includegraphics[scale=.5,trim= 10 20 0 20,clip]{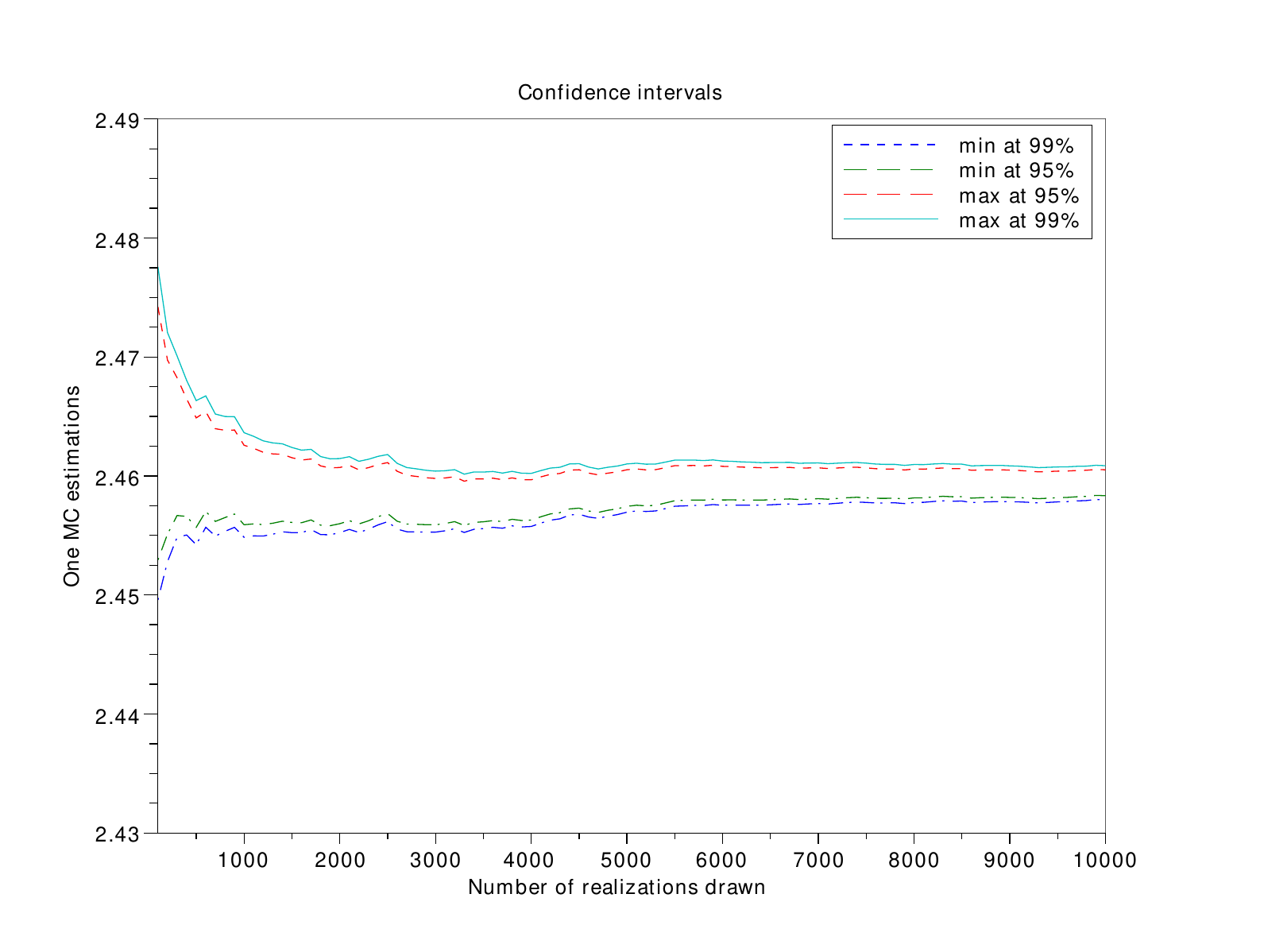}}
%   \end{tabular}
  \end{tabular}
\caption{\label{fig:MCconfidence} 
Top: maximum, mean and minimum in a sample of MC estimations $V_{M_{\rm test}}(Z^\lambda-\sum_{i=1}^I\alpha_i^\lambda\tilde Y^i)$ at various $\lambda$, as a function of $I$. 
The decrease rate is less fast among parameter values $\lambda$ not used by the greedy algorithm, but it is still reasonably good.
Middle: confidence intervals with probability 99\% and 95\% for the MC estimation $E_{M_{\rm test}}(Z^\lambda-\sum_{i=1}^I\alpha_i^\lambda\tilde Y^i)$ of $E(Z^\lambda)$ as a function of $I$ at one $\lambda$.
Bottom: confidence intervals with probability 99\% and 95\% for the MC estimation $E_{M}(Z^\lambda)$ of $E(Z^\lambda)$ as a function of $M$ at same $\lambda$, with a much slower decrease rate than above at a similar computational cost.}
\end{figure}

Moreover, our algorithm is robust with respect to the specific choice of $\tilde\Lambda$ above. That is, the control variates identified previously can still be used to reduce the variance of $E_M(Z^\lambda)$ at other $\lambda$ in $\Lambda$ than those in the trial sample $\tilde\Lambda$ used by the greedy algorithm above. Let us evaluate the gain for real-time purposes or in the limit of many queries in $\lambda$ when one can forget all the precomputations including the greedy construction. Then the computational cost per $\lambda$ is approximately divided by $10^4/10=10^3$, provided the reduced variance is still approximately of magnitude $10^{-4}$ for all $\lambda\in\Lambda$. This observation even holds if we require a higher precision (hence smaller tolerances $TOL$, $tol$ and correponsdingly increase  ``offline'' the dimension $N$ of the RB space for the PDE surrogates). Then, the gain grows like $TOL^{-2}$ in the infinitely-many-query asymptotics as long as one can reduce the variance by adding finitely-many control variates; and like $TOL^{-2}(1-I\times TOL^{-2}/Card(\tilde\Lambda))$ if one takes into account the greedy construction for a finite sample of $Card(\tilde\Lambda)\gg I\times TOL^{-2}$ parameter values, where $I$ is the minimal number of control variates required to achieve a reduced variance of magnitude $TOL$ in $\tilde\Lambda$. (Note that the computational cost of one realization becomes even more dominant as $N$ increases.) Here, Fig.~\ref{fig:MCconfidence} shows that that our greedy construction is robust: with control variates built using the $10\times 10$ trial values of $\lambda$ above, we did thorough MC estimations of the variance for $100$ other parameter values. In our example, it roughly holds $I\le-\log(TOL)$ for $TOL> 10^{-7}$ and the computational gain grows exponentially fast like $TOL^{-2}(1+TOL^{-2}\log(TOL)/Card(\tilde\Lambda))$ with the number $-\log(TOL)$ of accurate digits required for sufficiently large $Card(\tilde\Lambda)$. %x=(1:7)';y=1:5;plot(x,2*x*ones(y)+log10(1+(10.^(2*x)).*log(10.^(-x))*(10^(-y))))

In any case, if one cannot reduce the variance (or equivalently, if $I$ needs to be huge for any useful variance reduction), there is -- almost -- no loss of effort in using the RB control-variate MC method compared with the naive MC method, since would simply observe the absence of variance reduction using $M_{\rm test}$ realizations and would next need to increase the number of realizations to the same number of realizations as the naive MC method. Only the computation of coefficients $\alpha_i^\lambda$ would have been unnecessary, which is not much compared with a large number of realizations.

Of course, the efficiency of our RB control-variate MC method is limited by the cost and accuracy of single realizations, like any MC method. This has nothing to do with its ability at reducing the variance (and thus controlling the number of realizations needed for accurate MC estimations), but it is a practical limitation everyone doing simulations is confronted with. For UQ in particular, if we use a smaller correlation length $\delta=0.05$ but require the same tolerance $TOL$, the ``variational crime'' analysis above requires a KL truncation at a much higher-order $K=70$. Now, it is practically impossible today (in 2012) to carefully inspect all $2+K=72$ directions with a fine grid of trial values for the ``offline'' construction of RB surrogates\footnote{
  The computations presented in this work have been done in a very reasonable time, from a few minutes to one hour on a single processor unit, but inspecting $72$ directions with $2$ points per directions, that is $2^72\approx 10^{21}$ points, would require yers of computations.
}. So, even though the RB control-variate MC method still works perfectly well with inaccurate RB surrogates $u_{\mathcal{N},K,N}$, there is little computational gain possible on the whole simply by reducing the variance when the statistical error is dominated by the KL truncation error (there is no point in getting accurate MC confidence intervals if the RB surrogates are not good enough). The RB control-variate MC method is not a definitive cure to the ``high-dimensions'' problem for PDEs with stochastic coefficients which we already mentionned previously.

All existing numerical approaches to UQ problems that invoke a KL decomposition of the input noise are limited in practice by the ``curse of dimensionality'' in $K$ anyway, in so far as they require numerical approximations of realizations of the PDE solution in a space of high dimension $K$. (Most often in UQ people only use $K\le10$ in practice~\cite{frauenfelder-schwab-todor-2005,doostan-ghanem-redhorse-2007,nouy-2007}.) So is our approach. But in any case, when one can tackle the problem of simulating high-dimensional realizations, say only by expensive precomputations, then using our RB control-variate MC method and the PDE solution as a black-box still makes sense: the simulation of one realization of $Z^\lambda$ at one given $\lambda$ is typically all the more expensive as the parameter size increases, for instance just by assembling the matrix of the discrete BVP here, and the gain obtained by variance reduction (if any) is in fact all the more interesting. Moreover, there exist other ways of simulating the law of $b$ that do not require a KL representation (thus a limitation in $K$), even extend to other random fields than~\eqref{KL}, and can be used with our RB control-variate MC method~\cite{efendiev-hou-luo-2006}. This might be a cure to some problems but we keep its investigation for future works of ours. Last, for the model problem here, we are in fact quite fortunate, as already noted in~\cite{boyaval-lebris-maday-nguyen-patera-2009}, that the KL spectrum has a fast decay. The first KL modes span all the solution space when the required precision is not too small. Then we can still build ``offline'' good enough surrogates $u_{\mathcal{N},K,N}$ with a very small training set of trial parameter values, as we are going to see in the next section. Furthermore, since the previous observation holds only for specific choices of the spatial correlations, before shifting to the numerical proof of efficiency for the RB control-variate MC method in another UQ framework, we are also going to mention one promising way of tackling the simulation problem with high-order KL truncations which is currently under study, and where the RB control-variate MC method can still be useful. Our conclusion of the first example is thus mostly concerned with an efficient combination of the standard RB method with the RB control-variate MC method in a UQ framework.

\subsection{RB control-variate MC with RB surrogates}
\label{sec:certifiedRBmethod}

We briefly recall the ``standard'' RB method used above (see also~\cite{Machiels2000,Prud'homme2002a,Rozza2008} e.g.) before discussing its specific use in UQ with a MC method like our RB control-variate MC method. The reader already expert in the standard RB method may want to skip the next subsection.

\subsubsection{The standard RB method}

The goal of the RB method is to reduce the computational cost of solutions $u(\mu)$ to a PDE parametrized by $\mu$ when $u(\mu)$ has to be computed many times for many values of $\mu$. Here we set $\mu=(k_2,\bar E,\Upsilon\sqrt{\lambda_1}Z_1(\omega),\ldots,\Upsilon\sqrt{\lambda_K}Z_K(\omega))$ and given a Hilbert space $X = H^1(\D)$ with inner-product
$$
(u,v)_X = \int_{\D} \nabla u \cdot \nabla v + \int_{\Gamma_{\rm B}} u v\,,\ \forall u,v \in X
$$
% $$
% (u,v)_X = \int_{\D} u\:v + \int_{\D} \nabla u \cdot \nabla v\,,\ \forall u,v \in X \,,
% $$
and with energy norm $ \|u\|_X = \sqrt{(u,u)_X} $, we denote
$$
a(u,v;\mu) = k_1 \int_{\D_1} \nabla u \cdot \nabla v + k_2 \int_{\D_2} \nabla u \cdot \nabla v 
+ \int_{\Gamma_{\rm B}} b \: u \: v %\,,\ \forall u,v \in X \,,
$$
a bilinear form %$a(\cdot,\cdot;\mu)$ 
that is elliptic $ \forall u,v \in X $. %on $X$.
Then, for parameter values $\mu$ in a given range, we consider the solutions $u(\mu)\in X$ to~\eqref{eq:weak}, such that
\begin{equation}
\label{eq:weak_general}
 a(u(\mu),v;\mu) = l(v) \,,\ \forall v \in X \,,
\end{equation}
with outputs $s(\mu)=l(u(\mu))$. In practice, we assume that one can compute a good approximation of $u(\mu)$ in a linear subspace $X_{\cal N}\subset X$ with large dimension $\mathcal{N}\gg1$ for any $\mu$. Moreover, here we discretize $b\approx b_K$ so the fully computable approximations also involve a discretization parameter $K$ and is denoted $u_{\mathcal{N},K}(\mu)\in X_{\cal N}$. Computing $u_{\mathcal{N},K}(\mu)$ for many $\mu$ is costly and the goal of the RB method is to construct a linear subspace $X_N\subset X_{\cal N}$ with small dimension $N\ll\mathcal{N}$ such that, for all $\mu$ in the given range, some approximations $u_{\mathcal{N},K,N}(\mu)\in X_N$ to $u(\mu)\in X$ can be computed faster. The RB method invokes linear approximations $u_{\mathcal{N},K,N}(\mu):= \sum_{n=1}^N \gamma_n(\mu) u(\mu_n)$ where the coefficients $\gamma_n(\mu)$, $n=1\ldots N$, are computed fast at any $\mu\in\Lambda$ by the Galerkin method in $X_N$: $m = 1\ldots N$
\begin{equation}
\label{eq:weak_generalRB}
\sum_{n=1}^N \gamma_n(\mu) a(u(\mu_n),u(\mu_m);\mu) = l(u(\mu_m)) % \,,\ \forall j = 1,\ldots,N 
.
\end{equation}

The Galerkin approximation error $\|u_{\mathcal{N},K}(\mu)-u_{\mathcal{N},K,N}(\mu)\|_X$ can be estimated {\it a posteriori} 
by a fully computable error estimator $\Delta_N(\mu)$
% \begin{multline}
% \label{eq:residualestimator}
% \|u_n(\mu)-u(\mu)\|_X \le \Delta_n(\mu) \\ = \frac{1}{\alpha_{LB}} \sup_{v\in X/\{0\}} \frac{ |a(u_{n}(\mu),v;\mu)-l(v)| }{ \|v\|_{X} }.,
% \end{multline}
\beq
\label{eq:estimator}
\|u_{\mathcal{N},K}-u_{\mathcal{N},K,N}\|_X \le \underset{\|v\|_X=1}{\sup}\frac{a(u_{\mathcal{N},K,N}(\mu),v;\mu)-l(v)}{\alpha_{LB}(\mu)}
\eeq
using a lower bound $\alpha_{LB}(\mu)$ for the coercivity constant of the bilinear form, which is useful for the certifiability of the reduction method as well as for the selection of adequate parameter values $\mu_n$, $n=1\ldots N$, in the range of interest. Here, given the choice of inner-product on $X=H^1(\D)$,  
% $\int_{\D} \nabla u \cdot \nabla v + \int_{\Gamma_{\rm B}} u \: v$ for $u,v\in H^1(\D)$,
one can compute explicitly and at little expense quite a sharp lower bound $\alpha_{LB}(\mu)=\min(k_1,k_2,\inf b)$ that is uniformly good {\sl with respect to the parameter values} and for any admissible parameter range.

A good {\it reduced basis} $u(\mu_1),\ldots,u(\mu_N)$ spanning $X_N$ can be found in practice with a greedy algorithm like~\eqref{eq:greedy}~\cite{buffa-maday-patera-prudhomme-turinici-2009,binev-cohen-dahmen-devore-petrova-wojtaszczyk-2010}. The standard procedure selects $N$ parameter values $\mu_1,\ldots,\mu_N$ {\it incrementally} in a trial sample $\Theta$ for $\mu$ as follows: given any $\mu_1$ we next define for $n=1,\ldots,N-1$
\beq
\label{eq:optimpb}
\mu_{n+1}\in\underset{\mu\in{\Theta}}{\rm argsup}\: \Delta_n(\mu) \,.
\eeq

Note that the bilinear form $a$ is affine with respect to functions of the parameters (viz. parametrized by $\mu$ only through coefficients in a linear combination of bilinear forms, recall~\eqref{eq:weak_KL}). This is useful to efficient RB implementations in practice since the parameter-independent matrices can be precomputed. It allows fast evaluations\footnote{
  Denoting by $M$ the matrix of the inner-product $(\cdot,\cdot)_X$ computed for the discrete FE space where the parametrized BVP reads $A(\mu)U(\mu)=B$, the RB approximation error in the ``energy'' norm $\sqrt{(U(\mu)-U_n(\mu))^TM(U(\mu)-U_n(\mu))}$ is evaluated by the residual estimator $\sqrt{(A(\mu)U_n(\mu)-B)^TM^{-1}(A(\mu)U_n(\mu)-B)}/\alpha_{LB}(\mu)$.
%   where $\alpha$ is a numerical approximation of the lowest eigenvalue of a generalized symmetric problem $A(\mu_0)V = \alpha MV$ computed 
%   (using standard Arpack routines for instance) at some $\mu_0$ well chosen {\it a priori} to yield a lower bound of $\alpha_{LB,\mathcal{N}}(\mu)$ at all $\mu$.
  Now, $(A(\mu)U_n(\mu)-B)^TM^{-1}(A(\mu)U_n(\mu)-B)$ is a quadratic form in the small-dimensional vector $U_n(\mu)$ containing the coefficients of the approximation in the reduced basis. And its entries are themselves quadratic in (functions of) $\mu$ when the entries of $A(\mu)$ depend linearly in (functions of) $\mu$. The latter can thus be precomputed for a given reduced basis  and next be assembled in $O(N^2 K^2)$ to yield fast error estimations at any $\mu$.
}  of $u_{\mathcal{N},K,N}(\mu)$ and $\Delta_N(\mu)$ at any $\mu$.
Here, the quality of our RB approximations is similar to that in~\cite{boyaval-lebris-maday-nguyen-patera-2009}.

\subsubsection{Certified RB control-variate MC}

First, note that to get certified results in the frame of MC simulations, one should check the quality of those RB approximations when they are used, since in any case the precomputed basis for the RB surrogates invoked at each $\lambda$ and for each realization has been trained ``offline'' only with a few trial parameter values in $\Theta$ (possibly forgotten afterwards). In the RB control-variate MC method, it is necessary to check the RB approximation error with~\eqref{eq:estimator} at least for the $M_{\rm large}$ realizations in $E_{M_{\rm large}}(Z^{\lambda_i})$, $i=1\ldots I$, at the parameter values $\lambda_i$ where control variates are constructed
plus, for each $\lambda$, at the $(1+I)M_{\rm test}$ realizations of $Z^{\lambda}$ and $Z^{\lambda_i}$, $i=1\ldots I$, that are used to compute the final MC estimation with reduced variance. In practice here, this increase of computational cost per realization and $\lambda$ % at least in the case where $K$ is not too large,
was not much (that is, only a fraction of the total computational time without certification).

Moreover, for $\delta=0.5$ and $K=10$, we precomputed a reduced basis of dimension $N=12$ for the PDE solutions using only a cartesian grid of $10\times10\times2^K$ values for $\mu$.
Although this is a coarse training sample, we never had to enrich the RB approximation space $X_N$ ``online'', that is during the MC simulations. Thus, for UQ problems similar to our example, although our MC approach does not aim at exploiting optimally the regularity of the solution, in particular because the choice of trial parameter values for a greedy algorithm is ``blind'', the RB control-variate MC method is clearly interesting: it is simple to implement, without much a priori knowledge of the solution, and accurate at a reasonable computational cost (``offline'' and ``online'' computations take a few minutes on a single processor unit here).

\begin{figure}
\centering
%  \begin{tabular}{r}
    \makebox[.5\textwidth][c]{\includegraphics[scale=.5,trim= 10 20 0 20,clip]{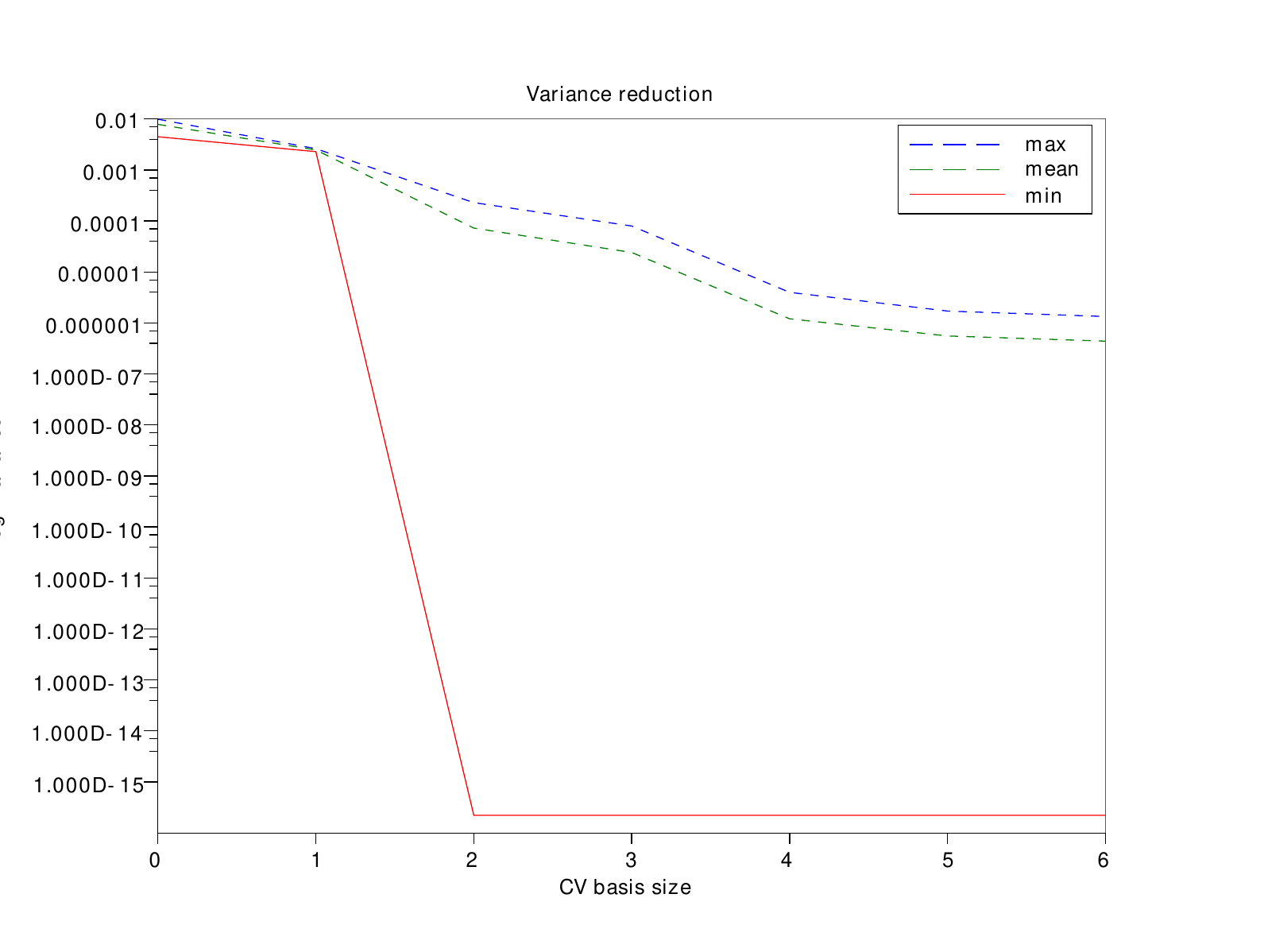}}
%  \end{tabular}
\caption{\label{fig:MC2} 
Maximum, mean and minimum in a sample of MC estimations $V_{M_{\rm test}}(Z^\lambda-\sum_{i=1}^I\alpha_i^\lambda\tilde Y^i)$ at various $\lambda$, as a function of $I$. 
The decrease rate is less fast among parameter values $\lambda$ not used by the greedy algorithm, and less fast than in~Fig.~\ref{fig:MCconfidence}.}
\end{figure}

For $\delta=0.05$, we mentionned previously that the problem is not to reduce the variance of MC estimations, but to get certified MC estimations when one uses a KL representation of the input noise at each of the MC realizations, because the parametrization of the latter needs approximately $K=70$ dimensions, which cannot be explored finely ``offline'' in a reasonable time (one manifestation of the well-known ``curse of dimensionality''). One way out is to simulate differently the noise, in a way that could invoke RB surrogates in a low-dimensional parameter space only~\cite{efendiev-hou-luo-2006}. We shall study this approach in future works. Another way if we want to use KL representations is, first, to observe that since the error per realization is averaged in $E(s_{\mathcal{N},K,N}(\lambda))$ one can still afford some realizations with large error as shown in~\cite{boyaval-lebris-maday-nguyen-patera-2009}. Moreover, since the solution space very much depends on the first KL modes only as observed in~\cite{boyaval-lebris-maday-nguyen-patera-2009}, one can build here ``offline'' surrogates $u_{\mathcal{N},K,N}$ in a space of low dimension. Here, we could use RB surrogates from a small linear space $X_N$ of dimension $N=20$ without having to enrich $X_N$ ``online'', although $X_N$ was built ``offline'' with a greedy algorithm exploring only a very coarse training sample of trial parameter values (using 10 trial values for each parameter component $\lambda_1=k_2,\lambda_2=\bar E$ and 2 trial values for each parameter component $\Upsilon \sqrt{\lambda_k} Z_k$, $k=1\ldots 10$, while for $k=11\ldots 70$ the parameter components were fixed).
% \footnote{
%    We stopped the greedy algorithm when a posteriori error estimates then satisfied $\|u_{\mathcal{N},K}-u_{\mathcal{N},K,N}\|_{H^1}\le\Delta_N\le.22\le\sqrt{.5}\le\sqrt{\alpha_{LB}}$ at all trial parameter values for the first time, which is a bit more careful and strigent than before.
% }. 
Then, compared with a direct ``naive'' MC simulation of the law of $Z^\lambda\equiv s_{\mathcal{N},K,N}(\lambda)$ at any $\lambda\in\Lambda$, the RB control-variate MC method provides computational reduction by reducing the variance of MC estimations like in the case where $\delta=0.5$, see Fig.~\ref{fig:MC2} in comparison with Fig.~\ref{fig:MCconfidence}. Note yet that the variance $V(Z^\lambda)$ is approximately one order of magnitude smaller when $\delta=0.05$ than in the previous example when $\delta=0.5$, and the variance decreases a bit more slowly than before, so the gain at same accuracy is a bit less than one order of magnitude smaller (and non zero only if we require a minimum precision on the confidence interval of the MC estimations of a bit more than one order of magnitude higher than when $\delta=0.05$).

% 
%  with the MC method or with the RB control-variate MC method and certify that the RB approximation errors due to the the RB surrogates in the MC estimations $E_{\rm M_{test}}(s_{\mathcal{N},K,N}(\lambda)-\hat Y^\lambda)$ with the control variate $\hat Y^\lambda$ for $Z^\lambda\equiv s_{\mathcal{N},K,N}(\lambda)$ is not too large. 
% Not in passing that 

Last, we would like to mention a way to reach higher KL truncation order $K$ which is not specific to our choice of the correlation length and where the RB control-variate MC method can still be useful. Note first that, quite often, surrogates $u_{\mathcal{N},K,N}(\lambda)$ built for a coarse training sample (say with only 1 trial parameter value in most directions) are already not too bad because the KL spectrum decays anyway. So, at any $\lambda$, one can expect $V(s_{\mathcal{N},K,N'}(\lambda)-s_{\mathcal{N},K,N}(\lambda))$ to be significantly smaller\footnote{
 In our numerical example above, $V_{\rm M_{test}'}(s_{\mathcal{N},K,N'}(\lambda)-s_{\mathcal{N},K,N}(\lambda))$ is in fact always close to zero machine by the property mentionned above:
at a moderately high precision level, the solution space is almost entirely spanned by variations in the first $K$ modes only. But of course this would not be true for random fields with thicker distribution tails.
} than $V(s_{\mathcal{N},K,N}(\lambda))$ and next improve the MC estimations of $E(s_{\mathcal{N},K,N}(\lambda))$ ``online''  at a very reasonable cost, by the combination of the RB control-variate MC method and the Multi-Level Monte-Carlo (MLMC) method~\cite{cliffe-giles-scheichl-teckentrup-2011} in two steps: estimate (i) $E(s_{\mathcal{N},K,N}(\lambda))$ with the RB control-variate MC method and (ii) $E(s_{\mathcal{N},K,N'}(\lambda)-s_{\mathcal{N},K,N}(\lambda))$ with a small number $M_{test}'$ of realizations. The dimension $N'\ge N$ of the enriched basis at $\lambda$ can be chosen a posteriori in (ii) from the few $M_{test}$ realizations in the MC estimation $E_{\rm M_{test}}(s_{\mathcal{N},K,N}(\lambda)-\hat Y^\lambda)$ of (i). In turn, the actual $I$ in (i) can be re-adjusted after (ii) to balance the statistical error in $E_{\rm M_{test}'}(s_{\mathcal{N},K,N'}(\lambda)-s_{\mathcal{N},K,N}(\lambda))$ such that $V_{\rm M_{test}}(s_{\mathcal{N},K,N}(\lambda)-\hat Y^\lambda)/{\rm M_{test}}\approx V_{\rm M_{test}'}(s_{\mathcal{N},K,N'}(\lambda)-s_{\mathcal{N},K,N}(\lambda))/{\rm M_{test}'}$ if possible. % how to choose $M_{test}',M_{test}$ ? = ?
Of course, the latter new strategy deserves a specific study, which we also keep for future works.

\section{Application to Bayesian estimation}
\label{sec3}

Let us now numerically demonstrate the efficiency of the RB control-variate MC method in another context useful for UQ, with parameterizations of higher dimension than in the previous section. As example, we consider the computation of a Minimum-Mean-Square-Error (MMSE) Bayes estimator in various parametric settings. To identify various possible cases with many queries in the parameter, we first use a toy-model. Second, we use the same model PDE problem~(\ref{eq:strong}--\ref{eq:BC}) as in the previous section. In the latter case, we also combine the RB control-variate MC method with standard RB surrogates of the PDE solutions, like in the previous section. But of course, the use of surrogate models in Bayesian estimation is not new~\cite{efendiev-hou-luo-2006,nguyen-rozza-huynh-patera-2009}.

\subsection{A toy model for various Bayesian estimations}
\label{sec:toy-bayesian}

The Bayesian estimation of some unknown input parameter $\theta$ in a model consists in improving a prior guess of its probability distribution $\pi^\xi$ (the probabilistic counterpart of the deterministic Tikhonov regularization in inverse problems, see~\cite{marzouk-najm-2009} e.g.) by observations of an output $s^\lambda(\theta)$ of the model that are viewed as realizations of a random variable~\cite{robert-2007}. We denote by $\lambda$ a control parameter entering the model, and by $\xi$ a so-called {\it hyper-parameter} in the Bayesian frame. The random variations of $s^\lambda(\theta)$ can be generated by {\it aleatoric noise} (for a given fixed value of $\theta$, the truth is stochastic and induces a distribution of value $s^\lambda(\theta)$) as well as the {\it epistemic} uncertainty about the parameter $\theta$ (our observations of the truth are noised).

When one knows the value of a control parameter $\lambda$ in a model, as well as $J$ observations $s_j^\lambda$, $j=1\ldots J$ ($J$ i.i.d. realizations of $s^\lambda(\theta)$) with likelihood $f^{\lambda,\zeta}(s_j^\lambda|\theta)$, the probability density of $s^\lambda$ knowing $\theta$), then Bayes formula allows one to compute the posterior distribution of $\theta$
\begin{equation}
\label{eq:posterior} 
% \pi^{\lambda,\zeta,\xi}(\theta|\{s_j^\lambda\}) = \pi^\xi(\theta)\prod_{j=1}^J 
% \frac{f^{\lambda,\zeta}(s_j^\lambda|\theta)}{\int f^{\lambda,\zeta}(s_j^\lambda|\theta)\pi^\xi(\theta) d\theta}
\pi^{\lambda,\zeta,\xi}(\theta|\{s_j^\lambda\}) \propto \pi^\xi(\theta)\prod_{j=1}^J{f^{\lambda,\zeta}(s_j^\lambda|\theta)}
\end{equation}
where $\zeta$ is another hyperparameter (entering the likelihood and not the prior like $\xi$). The posterior is used to compute the probability distribution of quantities that depend on $\theta$. It is also used to compute a deterministic approximation $\hat\theta^{\lambda,\zeta,\xi}(\{s_j^\lambda\})$ of $\theta$ given $\{s_j^\lambda,j=1\ldots J\}$, for instance the Minimum-Mean-Square-Error (MMSE) estimator~\cite[p.349]{casella-berger-1990}, which in decision theory\footnote{ 
  Other risk functions than the quadratic loss and other decision rules than the risk minimization are also used in practice~\cite{casella-berger-1990,robert-2007}, which lead to different estimators for $\theta$. For instance, one also uses the expected linear loss function $E(|\bar\theta-\theta|)$ which is minimal at the median $\bar\theta$, defined by $P(\theta\le\bar\theta|\{s_j^\lambda\})=\frac12$, and the maximum likelihood principle that takes the maximum a posteriori $\max_{\theta}P(\theta|\{s_j^\lambda\})$ as estimator. In this work, we consider only the MMSE estimator~\eqref{eq:MMSEestimator}. Notice that it is a good example for applications where one is interested by the expectation of a smooth functional of $\theta$ in the end, because its computational complexity with a MC method if of the same kind.
  % even though it is not always the best interesting one.
  % We intend to extend the reduced-basis paradigm to other estimators in future works, where other samplings than the control-variate MC method should also be considered.
} is interpreted as the minimum of the expected quadratic loss $E(|\hat\theta^{\lambda,\zeta,\xi}(\{s_j^\lambda\})-\theta|^2)$, averaged over both distributions of $\theta$ and $\{s_j^\lambda, j=1\ldots J\}$ 
\begin{equation}
\label{eq:MMSEestimator}
\hat\theta^{\lambda,\zeta,\xi}(\{s_j^\lambda\}) %= E(\theta|\{s_j^\lambda\}) % CONFUSING
:= \int\theta\:\pi^{\lambda,\zeta,\xi}\left(\theta|\{s_j^\lambda\}\right)d\theta\,.
\end{equation}

In practice, expectations using~\eqref{eq:posterior} like~\eqref{eq:MMSEestimator} are typically computed many times:
\begin{itemize}
 \item for many values of the control parameter $\lambda$ entering the model for the ``truth'',
 \item using many different sets of observations $\{s_j^{\lambda_0}, j=1\ldots J\}$, typically at $\lambda_0=\lambda$ but not necessarily, and
 \item for many values of the hyperparameters $\xi,\zeta$ entering the Bayesian frame
(prior and likelihood are often chosen in -- conjugate -- parametric families of distributions).
\end{itemize}
Think of the various pieces produced in a factory, or of real-time estimation procedures, where one may want to vary all these parameters~! We would thus like here to accelerate the computation of many parametrized MMSEs in many-query frameworks with our RB control-variate MC method and numerically explore the various parametrized settings with a toy-model first. 

\begin{figure}
  \centering
%\vspace{-4cm}
  \begin{tabular}{c}
 \makebox[.2\textwidth][c]{\includegraphics[scale=.35,angle=-90,trim= 25 25 0 25,clip]{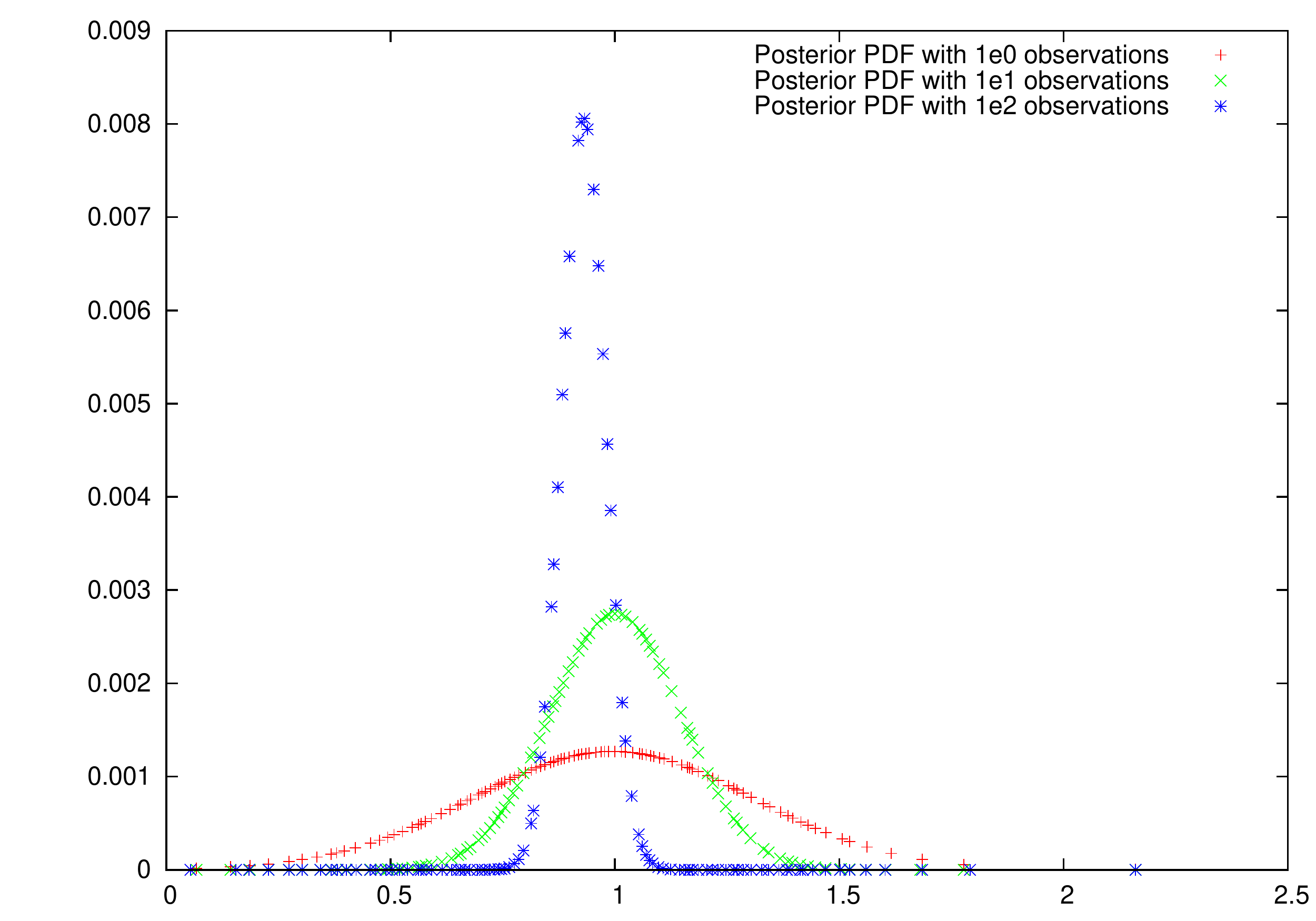}}
\\
%      \makebox[.4\textwidth][c]{\includegraphics[scale=.4]{fig/Elambda1.pdf}}
%    & \makebox[.4\textwidth][c]{\includegraphics[scale=.4]{fig/Vlambda1.pdf}}   
% \\
%      \makebox[.4\textwidth][c]{\includegraphics[scale=.4]{fig/Elambda0.pdf}}
%   &  \makebox[.4\textwidth][c]{\includegraphics[scale=.4]{fig/Vlambda0.pdf}}
% \\
%      \makebox[.4\textwidth][c]{\includegraphics[scale=.4]{fig/Elambda2.pdf}}
%    & \makebox[.4\textwidth][c]{\includegraphics[scale=.4]{fig/Vlambda2.pdf}}
% \\
     \makebox[.2\textwidth][c]{\includegraphics[scale=.45,trim= 30 05 50 40,clip]{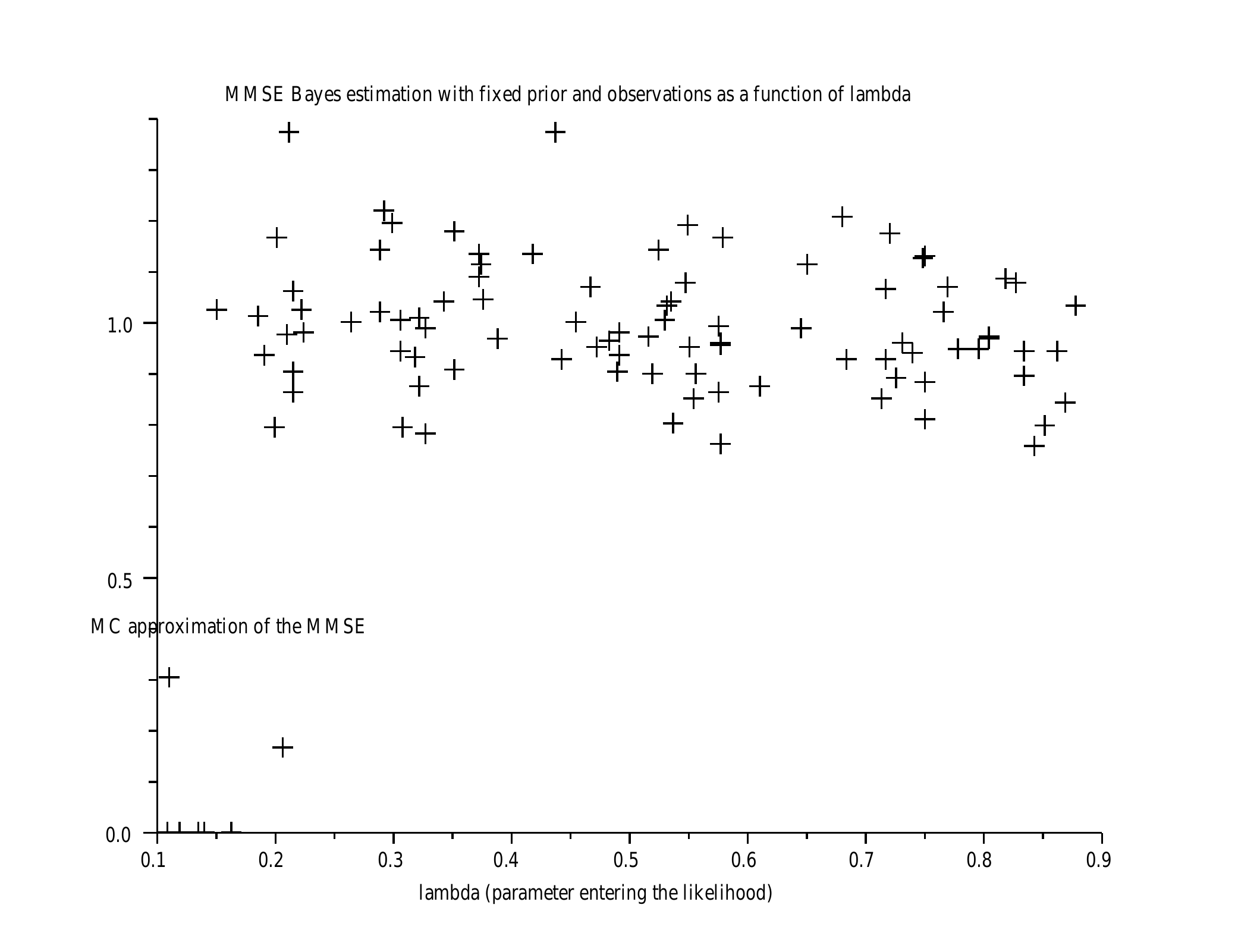}}
   \\ \makebox[.2\textwidth][c]{\includegraphics[scale=.45,trim= 30 25 50 40,clip]{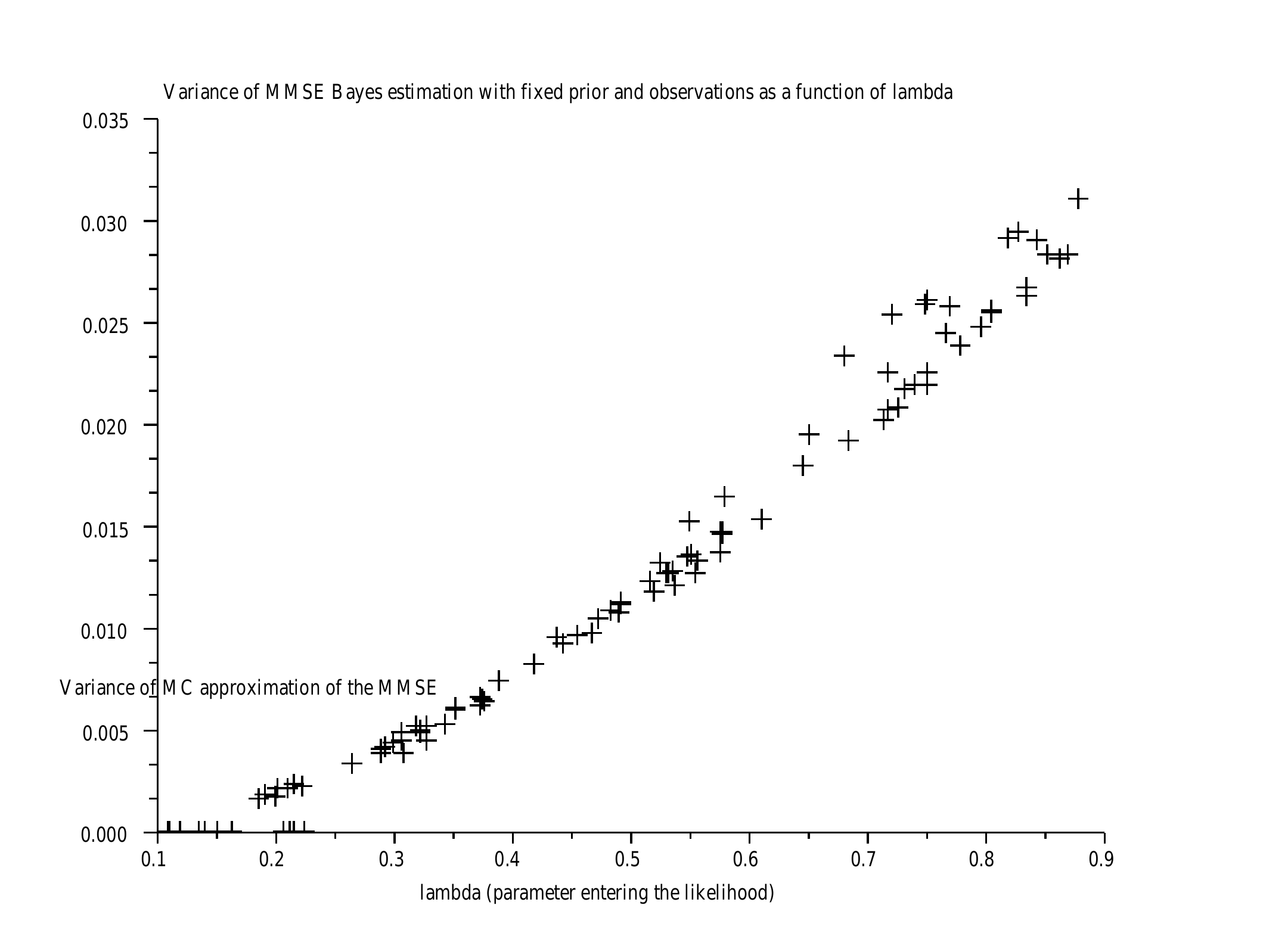}}
  \end{tabular}
\caption{\label{fig:lambda}
Posterior distributions computed for $J=1,10,100$ i.i.d observations $s_j\sim\mathcal{N}(\theta_0=1,\lambda^2=.25)$ (top) ;
and variations of the MC MMSE Bayes estimator (middle)
and its variance (bottom) as a function of $\lambda$ for $10^3$ realizations. %in~\ref{toy4}
%various parametrized scenarii: D.1, D.2, D.3 and D.4 from top to bottom. % A, B, C and D, 
%See the variance reduction in Fig.~\ref{fig:variance}.
}
\end{figure}

Let us choose a Bayesian frame for supposedly Gaussian observations of $\theta$ with known mean-square deviation $\lambda$, typically
$$
%s^\lambda(\theta)\sim\mathcal{N}(\theta,\lambda^2)
f^{\lambda,\zeta}(s|\theta) \propto e^{-\frac{|s-\theta|^2}{2\lambda^2}}
\qquad
\pi^{\xi=(\mu,\sigma^2)} = \mathcal{N}(\mu,\sigma^2) \,.
$$
Gaussian prior and likelihood are conjugate one-another %because the product of two Gaussian laws is still Gaussian
and the MMSE is analytically computable
$$
\hat \theta^{\lambda,\zeta,\xi}(\{s_j^\lambda\}) = \frac{\sigma^2}{\sigma^2+\lambda^2/J}\left(\frac1J\sum_{j=0}^J s_j^\lambda\right)
+ \frac{\lambda^2/J}{\sigma^2+\lambda^2/J}\mu
$$
with variance $\frac{\sigma^2\lambda^2/J}{\sigma^2+\lambda^2/J}$ (see~\cite[p.353]{casella-berger-1990} e.g.).
Note that for this simple model, there is no degree of freedom for a hyperparameter $\zeta$, but one already sees why one may want to use various values of the hyperparameter $\xi$. Although the MMSE Bayes estimation is asymptotically consistent and efficient (in fact, normally distributed) as $J\to\infty$, the quality of the estimation is clearly impacted by the choice of the hyperparameter at finite $J$. See the posteriors computed with $J=1,10,100$ observations and $\theta_0=1,\lambda=.5,\mu=.9,\sigma=.4$ in Fig.~\ref{fig:lambda}. %~\ref{fig:posterior}.

More generally, one can show under regularity assumptions on the likelihood $f^{\lambda,\zeta}(\cdot|\theta)$ (see~\cite[p.490]{lehmann-casella-1998} {\it e.g.}) that, if the observations $s^\lambda_j$ are indeed realizations of a random variable $s^\lambda(\theta)$ with density $f^{\lambda,\zeta}(\cdot|\theta)$ for one fixed value $\theta = \theta_0$ (hence distributed only due to the aleatoric noise of the model), then the MMSE %Bayes estimator 
converges in distribution %to the normal law 
when the number $J$ of observations increases
$$
\sqrt{J}(\hat\theta^{\lambda,\zeta,\xi}-\theta_0) \xrightarrow[]{d} \mathcal{N}(0,I(\theta_0)^{-1}) \,,
$$
with variance given by the {\it Fisher information}
$I(\theta_0) = E\left( \left( \partial_\theta \log f(\{s_j^\lambda\}|\theta) \right)^2 |\theta = \theta_0 \right) \,.$
The MMSE estimator is thus {\it asymptotically consistent} and {\it asymptotically efficient} and the choice of the prior becomes irrelevant as $J\to\infty$. Yet, for finite sample of observations, a good MMSE Bayes estimation strongly depends on the choice of the prior, in particular because MMSE Bayes estimators are biased. So one may try to optimize the hyperparameters given training samples of observations, possibly for various values of the control parameter: a natural many-query problem for the MC computation of parametrized MMSEs~! Let us see how the RB control-variate MC method performs in various parametrized settings.

\begin{figure}
  \centering
  \begin{tabular}{c}
    \makebox[.3\textwidth][c]{\includegraphics[scale=.3,trim= 0 50 0 50,clip]{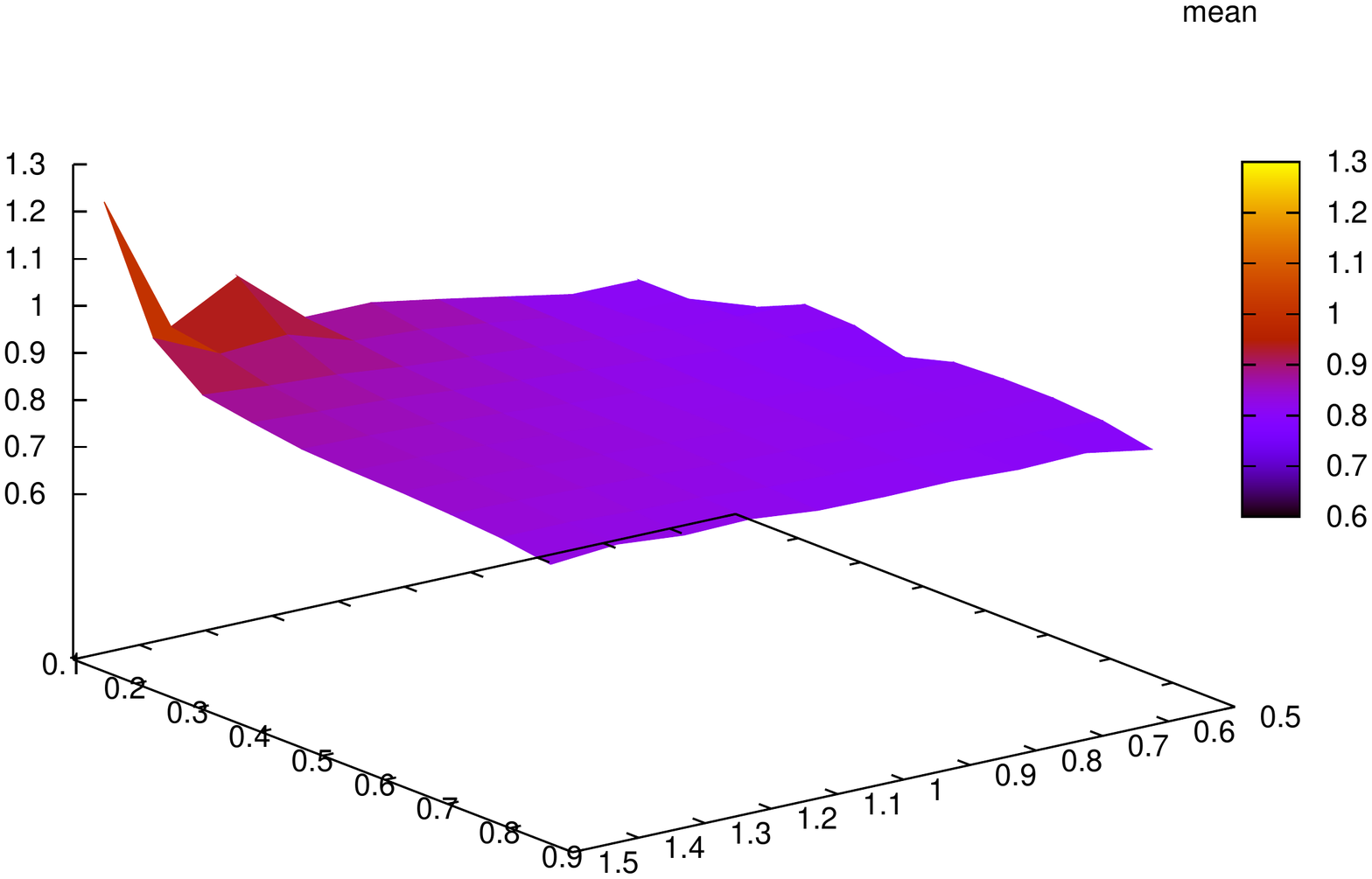}}
\\   \makebox[.3\textwidth][c]{\includegraphics[scale=.3,trim= 0 50 0 50,clip]{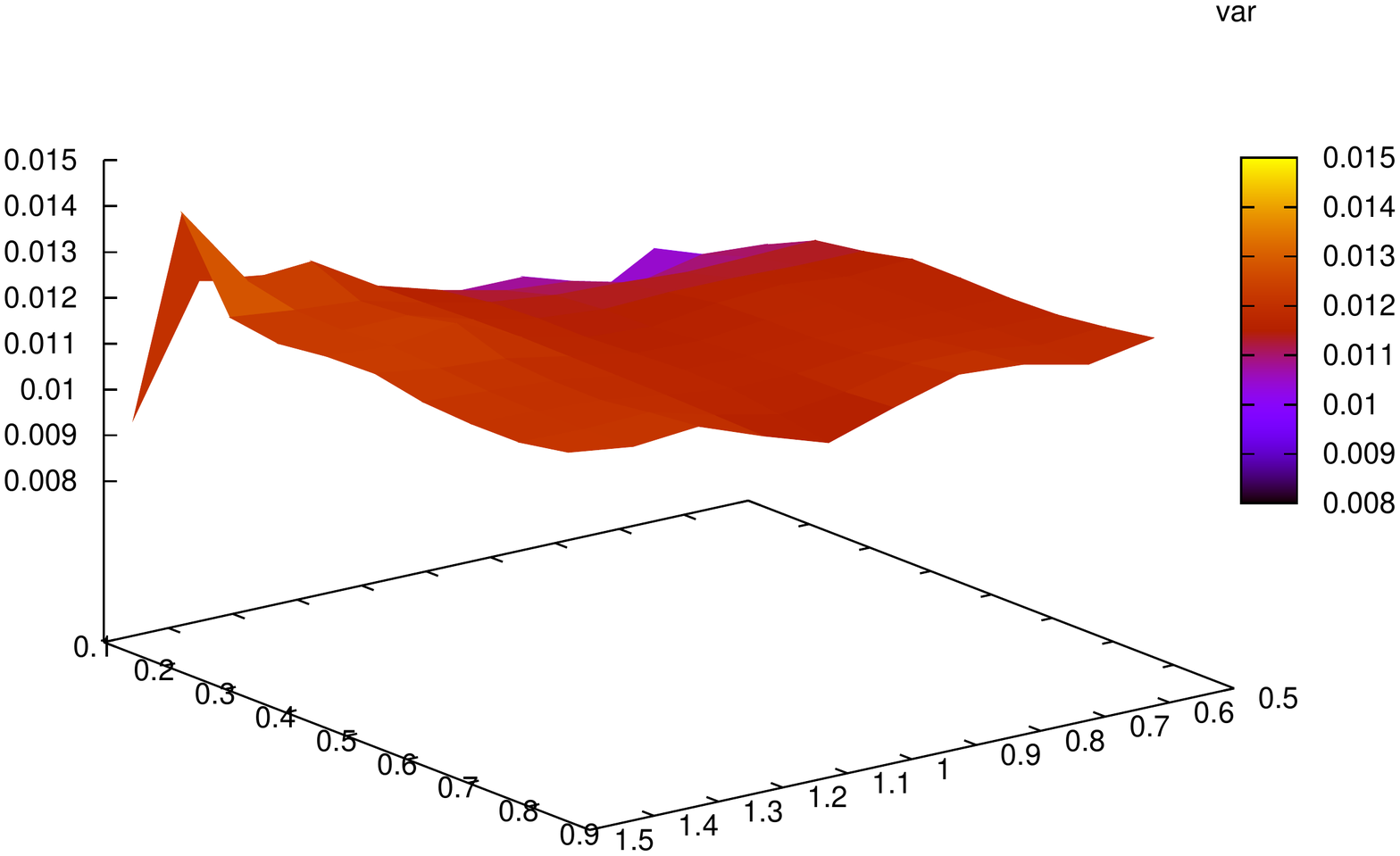}}   
\\ \makebox[.5\textwidth][c]{\includegraphics[scale=.4,trim= 10 10 10 10,clip]{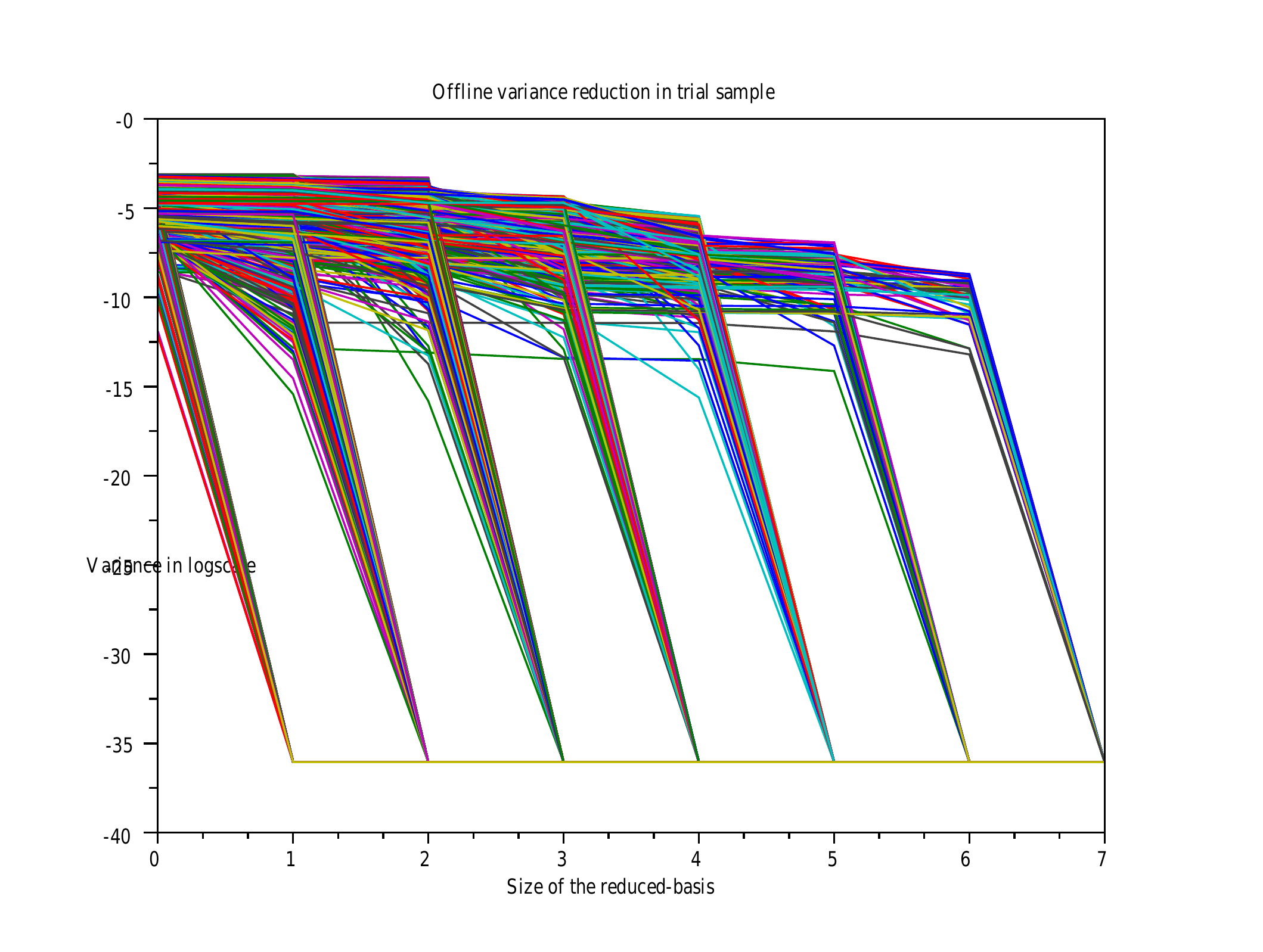}}
  \end{tabular}
\caption{\label{fig:delta}
Variations of MC MMSE Bayes estimation (top) and their non-reduced empirical variance (middle) %(with $10^3$ realizations) 
as functions of $(\mu,\sigma)$, using one set of observations and one $\lambda$. 
Bottom: empirical MC estimations of the reduced variances as functions of the number $I$ of variates for $20^3$ values of the parameter $(\{s_j^{\lambda_0},j=1\ldots J\},\lambda,\zeta\equiv(\mu,\sigma))$. Each colored line is the log-variance for one parameter value.}
\end{figure}

We do many MC estimation of the MMSE with $J=10$ i.i.d. synthetic observations generated with distribution law $\mathcal{N}(\theta_0,\lambda_0^2)$ and $\theta_0=1,\mu=.9,\sigma=.4$. We tested the RB control-variate MC method for various meaningful parametric variations:
\begin{itemize}
 \item $\lambda\in(.1,.9)$ with one set of observations generated at $\lambda_0=.5$ fixed or at $\lambda_0=\lambda$,
 \item $\lambda\in(.1,.9)$ and many sets of observations generated at $\lambda_0=.5$ fixed or at $\lambda_0=\lambda$,
 \item $\xi\equiv(\mu,\sigma)\in(.5,1.5)\times(.1,.9)$, $\lambda\in(.1,.9)$ and many sets $\{s_j^{\lambda_0},j=1\ldots J\}$ of observations generated at $\lambda_0=.5$ fixed or at $\lambda_0=\lambda$.
\end{itemize}
We always obtained meaningful estimations, and good variance reductions (more than $10$ orders of magnitude) with only a few variates (less than $10$). We show in Fig.~\ref{fig:delta} the ``most difficult'' case, which is still of course quite easy because of the very smooth dependence of our toy-model on the parameter. But this shows at least that even in potentially high-dimensional parameter contexts, there is some hope for the RB control-variate MC method to be useful.

\subsection{A PDE model with uncertain coefficients}

We now consider the PDE problem~(\ref{eq:strong}--\ref{eq:BC}) as a parametrized numerical model for Bayesian estimation. It is discretized like in the previous section: given $\delta$, one fixes the discretization parameters $\mathcal{N}$, $K$ and $N$. Note that we invoke the standard RB method~\cite{Prud'homme2002a,patera-rozza-2007} to compute fast numerous solutions $u$ to~(\ref{eq:strong}--\ref{eq:BC}) at each of the numerous values of the PDE coefficients invoked, which has already been used in a Bayesian estimation context with a deterministic quadrature formula for the output expectation~\cite{nguyen-rozza-huynh-patera-2009} rather than MC.

We still denote $\lambda=(k_2,\bar E)$ the parameter in the model. In addition to $s^\lambda$ defined in~\eqref{eq:output}, we also consider the temperature averaged at the top of the fin $o^\lambda=\int_{\Gamma_{N}\cap\overline{\D_1}}u$ as output. Clearly, there are many situations where the uncertain coefficient $b$ can a priori be modelled using only some rough representation as a random field~\eqref{KL}. Then, it often needs improving in a real setting using data from observations. This is the case of the permeability field for the Darcy equations in~\cite{efendiev-hou-luo-2006} e.g. 

Here, we shall try to improve i.i.d. {\it Gaussian} priors\footnote{
  For well-posedness of the PDE we will need to numerically truncate the Gaussian realizations by excluding those realizations where $b\le0$. This is a poor man's ``truncated Gaussian'', but in practice the standard deviations $\sigma_k$ are chosen sufficiently small for most $Z_k$ realizations to be in the range of admissibility where the surrogate model is valid.
} $Z_k\sim\pi^\xi\equiv\mathcal{N}(0,\xi_k)$, $\xi_k\equiv\sigma_k^2$, for the random variables $Z_k$, $k=1\ldots K$, using $J$ i.i.d. observations $o_j^{\lambda_0}$, $j=1\ldots J$, of the ouput $o^{\lambda_0}$ at some value $\lambda_0$ of the control parameter. Assuming one knows that the likelihood $f^{\lambda,\zeta}(o_j^{\lambda_0}|\{Z_k\})$ of the observations $o_j^{\lambda_0}$ given a realization $\{Z_k\}$ of the set of random variables $Z_k$, $k=1\ldots K$, are also Gaussian $f^{\lambda,\zeta}(o_j|\{Z_k\})\varpropto exp(-|o_j-o^\lambda(\{Z_k\})|^2/\zeta)$, $\zeta\equiv\sigma_O^2$, then the posterior distributions of the independent random variables $Z_k$ can be used to describe the specific context associated with the observations $o_j^{\lambda_0}$, $j=1\ldots J$. (They are a priori better than the prior guesses and should become all the better as the number of observations increases.) In particular, the model posterior, given by the Bayes formula $\prod_{k=1}^K \pi_k^{\lambda,\zeta,\xi}(Z_k|\{o_j^{\lambda_0}\}) \varpropto \prod_{j=1}^J f^{\lambda,\zeta}(o_j^{\lambda_0}|\{Z_k\}) \prod_{k=1}^K \pi_k^\xi(Z_k)$, can be used to compute various quantities of interest, like the expectation of $s^\lambda$. Let us use the RB control-variate MC method to evaluate fast $E_{\pi^{\lambda,\zeta,\xi}}(s^\lambda|\{o_j^{\lambda_0}\})$ for various values of the hyper parameters $\xi,\zeta$, various values of the control parameter $\lambda$ in the model and possibly various sets of observations (at $\lambda_0=\lambda$ possibly, but not necessarily). 
% Furthermore, one may want to reiterate the Bayesian estimation with the parametrized posterior $\prod_{j=1}^J f^{\lambda,\zeta}(o_j^{\lambda_0}|\{Z_k\}) \prod_{k=1}^K \pi^\xi(Z_k)$ when the control parameter $\lambda=(k_2,\bar E)$ in the numerical model varies as well as the set of observations $\{o_j^{\lambda_0}\}$, where possibly $\lambda_0=\lambda$.
Quite importantly, note that here we assumed that we know explicitly the likelihood function $f^{\lambda,\zeta}(\cdot|\{Z_k\})$ by experience. This is a specific Bayesian framework which is often used, although the exact likelihood function is often not known in practice (see Rem.~\ref{rem:precompute}). 

% For $\delta=.05$, $K=1$, $N=5$, let us vary the parameters as follows $\left(k_2,\bar E,\sigma_Z^2,\sigma_0^2 \right) \in [.1,10]\times[.1,1]\times[.03,.12]\times\{.02\}$ while observations $o_j^{\lambda_0}$,$j=1\ldots J$, are generated synthetically for $J=10$ by uncertainty propagation from $Z_1$ to $o^{\lambda_0}$ at $\lambda=\lambda_0$. The RB control-variate MC method used with a sample of $2\times 2 \times 2 \times 1\times 2^{10} = 8192$ parameter values
% % Each control variate is built using the same $M_{\rm large}=10^3$ MC realizations,
% % while a fast estimation of the variance is computed with $M_{\rm small}=10$ points during the greedy algorithm.
% allows to reduce the variance of naive MC estimations by a factor of $10^{10}$ with only $3$ control variates as shown in Fig.~\ref{fig:variancebayes1}, which shows the interest of the RB control-variate method in a high-dimensional context. Even though the model here is very simple, it applies in a Bayesian estimation context that is also used with more complex models, and it could thus be interesting then too.
% \begin{figure}
% \centering
% \makebox[.36\textwidth][c]{ \includegraphics[scale=.45,trim= 0 40 10 30,clip]{fig/variancebayes1.pdf} }
% \caption{\label{fig:variancebayes1} Maximum and mean of the reduced variances estimated by a fast MC during greedy inspection of the offline sample with cardinality $8192$ as a function of the number of control variates.}
% \end{figure}

\begin{figure}
\centering
  \begin{tabular}{c}
\makebox[.5\textwidth][c]{ \includegraphics[scale=.5,trim= 10 10 10 10,clip]{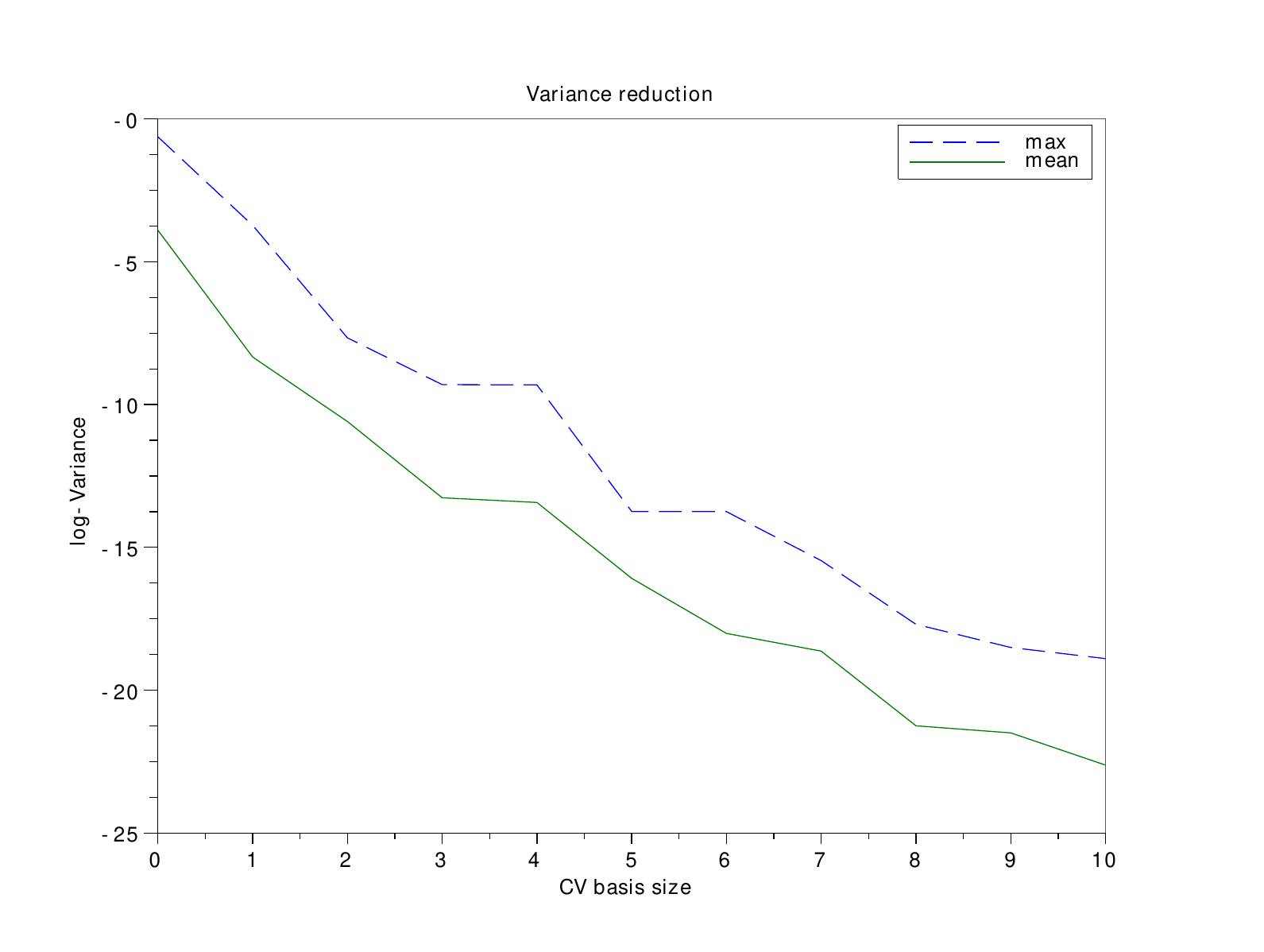} }   
\\
\makebox[.5\textwidth][c]{ \includegraphics[scale=.5,trim= 10 10 10 10,clip]{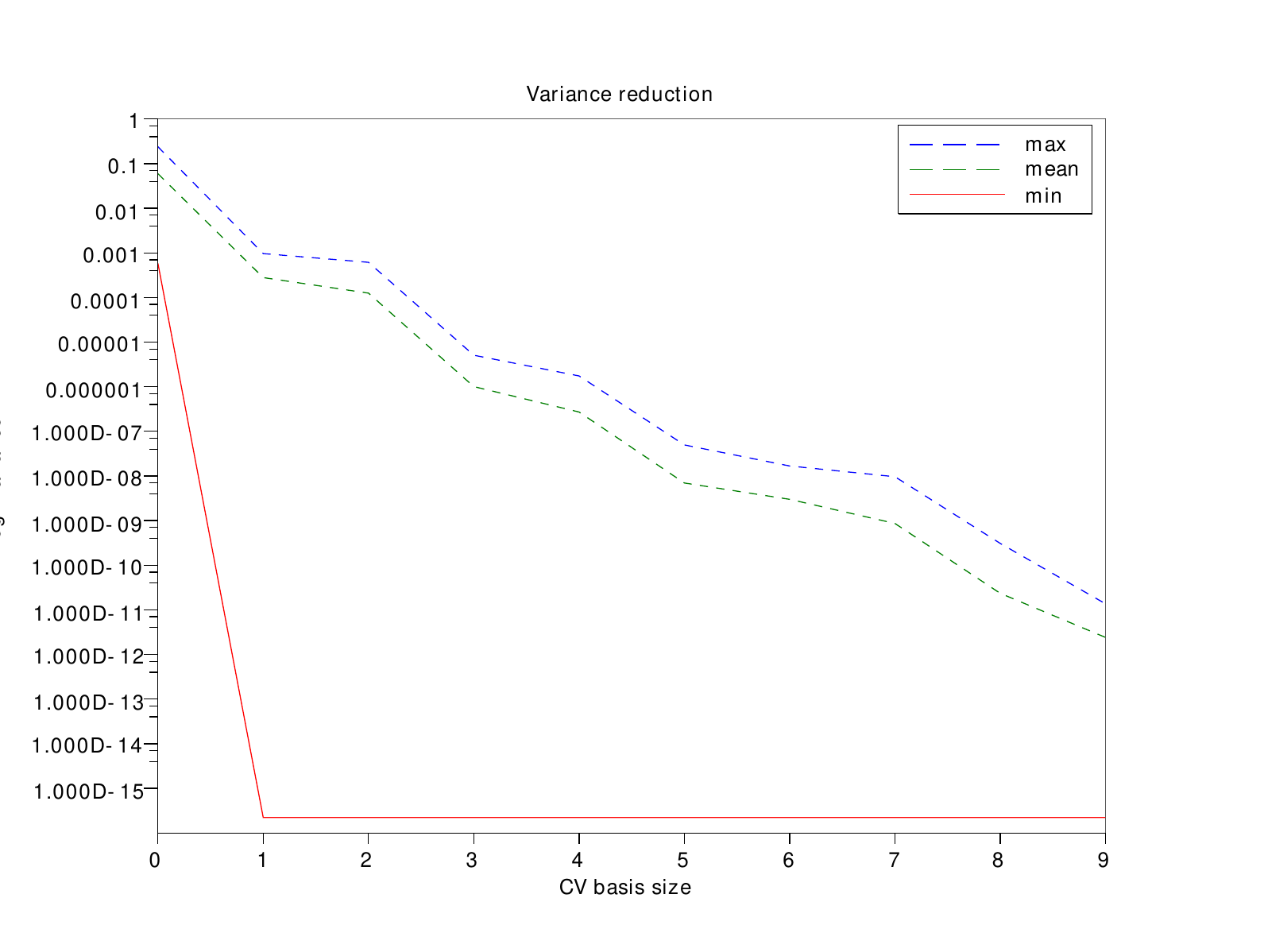} }
  \end{tabular}
\caption{\label{fig:meanVbayes1MCtest10} Top: maximum and mean in a sample of MC estimations $V_{\pi^{\lambda,\xi,\zeta}}(s^\lambda-\sum_{i=1}^I\alpha_i^\lambda\tilde Y^i|\{o_j^{\lambda_0}\})$ for various training values of $\lambda,\xi,\zeta$ and $\{o_j^{\lambda_0}\}$, as a function of $I$ during the greedy search (``offline'' stage) when $\delta=.5$, $K=10$, $N=12$ in case (ii) (case (i) is similar, with a decrease rate in fact a bit faster). Bottom: maximum, mean and minimum of the same estimators for various trial values of $\lambda,\xi,\zeta$ and $\{o_j^{\lambda_0}\}$ (average using many realizations) during the ``online'' stage.}
\end{figure}

We use the RB control-variate MC method to compute the MMSE of $s^\lambda$ with the parametrized posterior above in different settings for $\delta=.5$, $K=10$, $N=12$. For instance: (i)  with $4\times 4\times 2^K=2^12$ parameter values $\left(k_2,\bar E,\sigma_k^2\right) \in [.1,10]\times[.1,1]\times[10^{-4},10^{-3}]$ and one fixed $J=1$ observation $o_1^{\lambda_0}$ generated synthetically by uncertainty propagation at $\lambda_0=(2,0.5)$ or (ii) with $10\times10$ parameter values $\left(k_2,\bar E\right) \in [.1,10]\times[.1,1]$ and $10$ values for each of the $J=3$ observations $o_j^{\lambda_0}$,$j=1\ldots J$, generated synthetically by uncertainty propagation at $\lambda_0=(2,0.5)$. In both cases, the RB control-variate technique can bring computational reductions to a direct MC method as one can see by confronting the numerical results in Fig.~\ref{fig:meanVbayes1MCtest10} with the reasoning detailed in Section~\ref{sec2}. 

\begin{remark}\label{rem:precompute}
Sometimes we do not know how to choose the likelihood function and we need to precompute it numerically first. Let us illustrate this, in another Bayesian context for the sake of simplicity in the computations. For instance, assume that we want to esimate probabilities for $k_2$ from sufficiently many observations of the output $s^\lambda$ whose random fluctuations are generated by those of $b$. To estimate $k_2$ from $J$ i.i.d. observations $s_j^{\bar E_0}$, $j=1\ldots J$, with likelihood $f^{\bar E_0}(s_j^{\bar E_0}|k_2^0)$ depending only on the law of $b$ at the unknown value $k_2=k_2^0$, we could use a Bayesian approach with a prior distribution $k_2\sim \pi$. The posterior distribution %of $k_2$
\begin{equation}
\label{eq:bayesdistrib} 
\pi^{\bar E_0}(k_2|\{s_j^{\bar E_0}\}) \propto \pi(k_2) \prod_{j=1}^J f^{\bar E_0}(s_j^{\bar E_0}| k_2)
\end{equation}
allows one to compute quantities that depend on $k_2$ more accurately than with the prior distribution. (At least if the assumptions of the Bayesian model are satisfied by the ``truth'' and if the number of observations $J$ is large enough.) But first we need to precompute the likelihood $f^{\bar E_0}(\cdot| k_2)$.

A ``kernel'' approach for instance allows one to precompute the parametrized Probability Density Functional (PDF) from a sample of realizations at specific values of the controlled parameter $\bar E_0$. It amounts to give to each realization of $s^{\bar E_0}$ a weight.
%(Note again that it is possible to get a good sample of PDE solutions realizations because we use a surrogate built with the standard RB method in the previous Section~\ref{sec2}.) 
In Fig.~\ref{fig:variance4}, we show the approximate likelihood for observations $s_j^{\bar E_0}$ at one parameter value $\bar E_0$, after reconstruction of the PDF from a sample of realizations of $s^{(k_2,\bar E_0)}$ at various values of $k_2$, synthetically generated with a numerical MC simulation of our model and weighted by a truncated Gaussian weight with hyperparameters $\alpha_1,\alpha_2,M$
%adjusted to the values $s_m(k_2,\bar E)$
\beq
\label{eq:pdf}
%PDF^{\bar E_0,\alpha_1,\alpha_2,M}(s|k_2) \propto 
\sum_{m=1}^{M} 1_{|s-s_m^{(k_2,\bar E_0)}|<\alpha_2}e^{-\frac{|s-s_m^{(k_2,\bar E_0)}|^2}{2\alpha_1^2}} \,.
\eeq
% Notice that computing the PDF for each parameter value is time-consuming, 
% so we precompute it for some $k_2,\bar E$ and next interpolate the value at any other $k_2,\bar E$. 
% This may be inaccurate for high-dimensional $k_2,\bar E$ but sufficient for Bayesian estimation.

We numerically checked that the RB control-variate MC method was still efficient\footnote{
  Efficiency is meant with respect to variance reduction. Of course $k_2$ is only a scalar in the example here and one can argue that most often the MC method is not the best integration method then. % to compute the integral~\eqref{eq:MMSEestim}. 
  %In any case, it is very robust (in the case the integrand is really not smooth), and 
  Our example is nevertheless an illustration (admittedly simple) of more general situations with high-dimensional $k_2$ when one uses a robust MC-like method.
  %where one would be interested on a smooth scalar functional of $k_2$.
}
to compute many expectations of a scalar functional of the uncertain parameter $k_2$ using~\eqref{eq:pdf} as likelihood in the posterior~\eqref{eq:bayesdistrib}. % probability density $\prod_{j=1}^J f^{\bar E}(s_j^{\bar E_0}| k_2) \pi(k_2)$, we simply consider the expected value of $k_2$ itself, 
% In particular, one can compute deterministic approximations of $k_2$ using the MMSE estimator simply written
% \begin{equation}
% \label{eq:MMSEestim}
% \hat k_2 = E_\pi(k_2|\{s_j\}) \,.
% \end{equation}
% % (The MMSE interpretation in decision theory is the minimum of the quadratic loss $E(|\hat k_2 - k_2|^2)$ averaged over both distributions of $k_2$ and $\{s_j\}$, see~\ref{sec:toy-bayesian}.)
We used it for many values of the control parameter $\bar E_0$, of the set of observations $\{s_j^{\bar E_0}\}$ % (when $\bar E = \bar E_0$ or $\bar E \neq \bar E_0$), 
and of the hyper parameters $\alpha_1,\alpha_2,M$. But notice that actually using~\eqref{eq:pdf} in the posterior~\eqref{eq:bayesdistrib} is very time-consuming, especially with large $M$. Moreover, one needs to interpolate~\eqref{eq:pdf} at the values of $\bar E_0$ where it has not been precomputed. So the computation of probability densities, for instance with a kernel approach, remains a challenging numerical problem in Bayesian estimation. And the computational reductions proposed in the present work do not tackle the latter point. But note that clearly, to improve the computation of~\eqref{eq:pdf} at many parameter values $\bar E_0$ in particular, RB ideas could still be useful. Yet this is still another topic for another work.
% Then, given a uniform prior $k_2\sim \mathcal{U}([.1,10])$ and the likelihood~\eqref{eq:pdf}, we show in Fig.~\eqref{fig:variance4} the variance reduction obtained within a trial sample of $10^3$ control parameter values $\bar E\in [.1,1]$ for MC MMSE estimation using $J=10$ observations $s^{\bar E_0}$ such that $\bar E = \bar E_0$ or $\bar E \neq \bar E_0$. The synthetic observations are generated for a fixed $k_2=k_2^0$ with a simple PDE model ($\delta=.5$, $K=5$, $N=25$). % estimated with $M_{\rm small}=1e2$ realizations
% % plus some measurement error distributed as a truncated Gaussian error with same hyperparameters $\alpha_1,\alpha_2$ as in~\eqref{eq:pdf}. 
% % The case $\bar E \neq \bar E_0$ is more difficult as shown in Fig.~\ref{fig:variance4}. ??
\begin{figure}
\centering
  \begin{tabular}{c}
     \makebox[.3\textwidth][c]{  \includegraphics[scale=.25,trim= 50 50 50 50,clip]{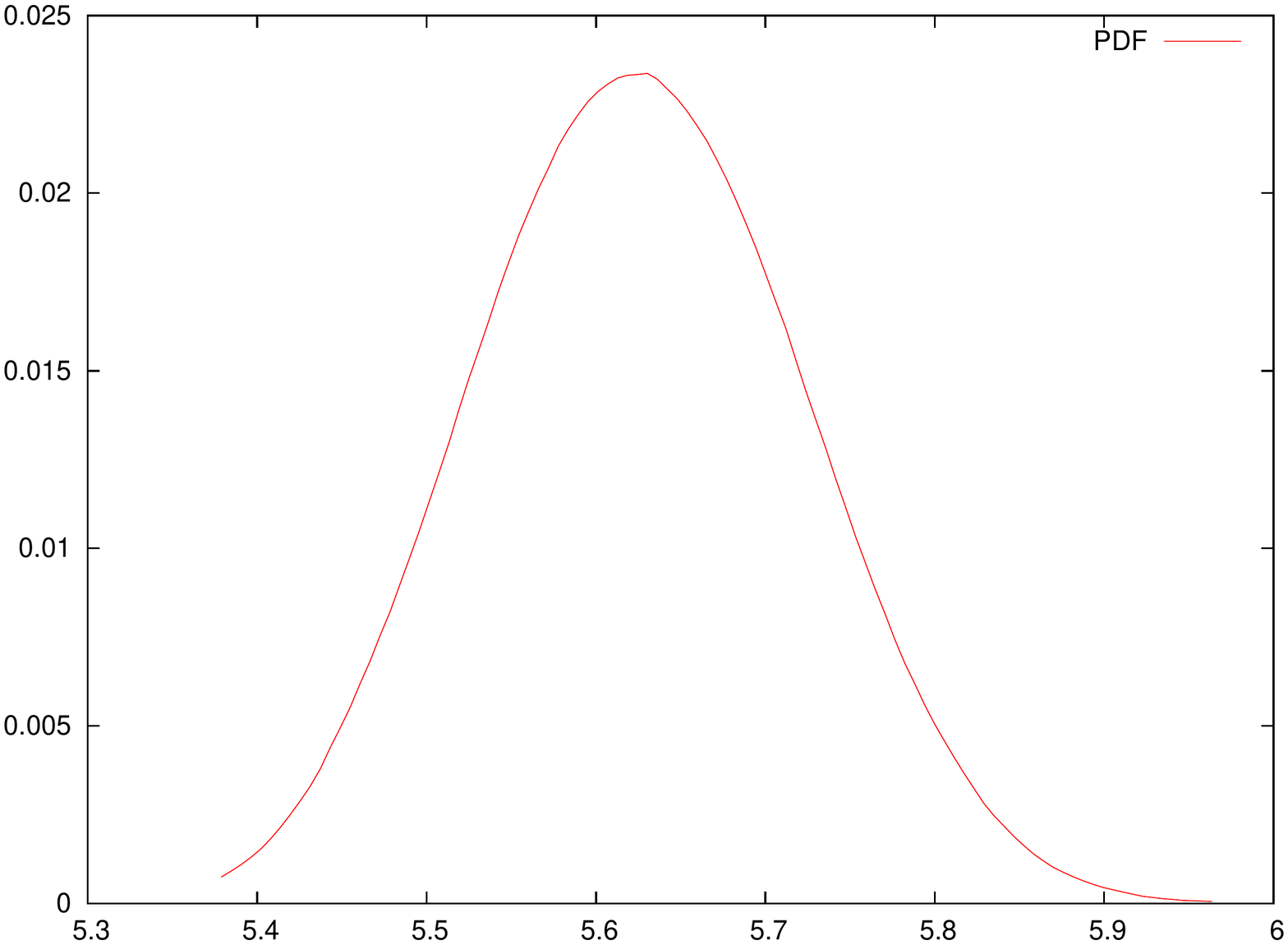} }
% \\    \makebox[.45\textwidth][c]{  \includegraphics[scale=.5,trim= 0 42 10 30,clip]{fig/variance4a.pdf} }
% %\\    \makebox[.45\textwidth][c]{  \includegraphics[scale=.5,trim= 0 42 10 30,clip]{fig/variance4b.pdf} }
% \\    \makebox[.45\textwidth][c]{  \includegraphics[scale=.5,trim= 0 42 10 30,clip]{fig/variance4c.pdf} }
  \end{tabular}
\caption{\label{fig:variance4} %From top to bottom:
Empirical PDF~\eqref{eq:pdf} of $k_2$ at fixed $\bar E_0$ for $M=10^4$ realizations with truncated Gaussian weight $\alpha_1=.5,\alpha_2=.01$.
% Reduced variances of MC MMSE estimators of $k_2$ using $J=10$ observations at $10^3$ parameter values $\bar E$ as a function of the number of control variates when $\bar E = \bar E_0$ ;
% similar results when $\bar E \neq \bar E_0$.
}
\end{figure}
\end{remark}

\section{Conclusion and Perspectives}

In this article, we have presented new elements of analysis of the RB control-variate MC method (error estimates and convergence results) as well as new useful applications (to uncertainty quantification of PDEs with stochastic coefficients). But a number of questions remains, about theory and practice, which may be the object of future works.

First, our convergence results are very much inspired from those in~\cite{buffa-maday-patera-prudhomme-turinici-2009,binev-cohen-dahmen-devore-petrova-wojtaszczyk-2010} about the standard RB method for PDEs and the same limitations apply: (i) they rely on assumptions that are difficult to check in practice (the fast decay of Kolmogorov widths), (ii) they may still be very pessimistic in comparison with practical possibilities (the constant in the upper-bound may be far suboptimal and a theoretical gain could be seen only for uninterestingly large dimensions $I$ of the reduced basis), and (iii) they do not help choosing a good training sample $\tilde\Lambda$ of trial parameter values for the greedy algorithm (especially in the cases where $\Lambda$ is large). Besides, (iii) is sometimes not only due to a lack of answer to the question how to het the optimal reduction, but also to the question how to get {\it any} reduction, when $\Lambda$ is large. Indeed, although the maximal gain (in the infinitely-many-query limit or, equivalently, for real-time purposes) reads only like the ratio of the reduced vs. non-reduced marginal computational time per parameter value $\lambda$, it is limited in practice by the possibility to inspect all trial parameter values $\lambda\in\tilde\Lambda$ and should take into account that ``offline'' part of the work. That is why, in absence of theory for choosing $\tilde\Lambda$ when $\Lambda$ has a high-cardinality (possibly infinite) {\it and} high dimension, one still needs numerical tests to check the capabilities of the reduction method using specific (heuristic) choices of $\tilde\Lambda$. Fortunately, we have been able to show numerically that some gain is possible in practical situations.

In some UQ frameworks, we could actually prove numerically that some gains were possible by applying the RB control-variate MC method. But as any numerical proof, this was achieved on one specific model problem, for which some limitations were also clearly seen. So one would of course still need characterizing precisely the limitations for another specific problem. More precisely, one would need to specify numerically the efficiency regime in a given UQ context where a parametrized expectation has to be computed many time for many values of the parameter, which depends on the number of queries in the parameter values, but also on the decrease rate of the variance with respect to the dimension of the reduced basis of control variates and on the required accuracy. Furthermore, this efficiency regime will be reachable only for those problems whose ``offline'' work is actually do-able (in a human lifetime) and is thus limited in the dimension of the parameter $\lambda$ by the size of $\tilde\Lambda$, as well as in the numerical complexity of the underlying model for random realizations by the computation time of one single realization. That latter limit also sets new challenges for future works in UQ: the simulation of high-dimensional noise (with high-dimensional KL truncation order $K$) in the case of PDEs with stochastic coefficients, and the fast computation of empirical likelihoods (say by kernel approach) in Bayesian estimation. We already mentionned some tracks for future works about these issues in the course of this work.

To conclude, this work is a direct improvement on~\cite{boyaval-lebris-maday-nguyen-patera-2009} which we believe to be a useful numerical approach to some UQ problems, in particular because it can bring huge computational gains at almost no cost to the simple and widely-used ``naive'' MC method.
% This is also a reasonably simple alternative to other numerical methods for uncertainty propagation that might theoretically better use the regularity of the PDE solution with respect to variations in the stochastic/controlled coefficients (when the problem is analytic like here) but are in practice often limited to small-dimensional random coefficients. 
Notice that in any case we can certify the approximation error in all steps of the computations (of course in a probabilistic sense as for the MC method) and there is thus no loss of rigor in our approach compared with the naive MC method. Moreover, we hope that our suggestion, after~\cite{nguyen-rozza-huynh-patera-2009}, to use RB ideas in such fields as Bayesian estimations, which are known to have computationally demanding applications, will bring some new thoughts in that community, who might for instance want to identify new ``many-query'' opportunities of computational reductions using greedy algorithms.

% Can this be systematically
% integrated to large computational platforms (SALOME e.g.). 
% We hope to tackle this challenge in a near future.
% Perspectives are:
% \begin{itemize}
%  \item on the methodological side, to use RB ideas within other computationally expensive
% discretization schemes in a parametrized framework,
% % for parametrized importance sampling (Girsanov theory) and MCMC, quantization
% to try improve RB-like parameter selection (hence the Greedy optimality) with an
% analysis of the variations in parameter space
% % inspired by approximation theory, analysis of variance -- ANOVA --, machine learning\dots
%  \item on the application side, to tackle higher-dimensional problems (Gaussian mixtures),
% time-dependent filtering problems (sequential Bayesian estimations), or
% the construction of probability distributions by Maximum Entropy methods ``a la Soize'' e.g.
% \end{itemize}

\appendix

\section{Weak greedy convergence}
\label{app:weakgreedy}

An important tool for the proof is Lemma~\ref{lem:cohen}, which is a straightforward variation of a result proved in~\cite[2.2]{binev-cohen-dahmen-devore-petrova-wojtaszczyk-2010}.
% only the line with ||gj|-|yj||\le|gj-yj| has to be changed
\begin{lemma}%(See \cite[2.2]{binev-cohen-dahmen-devore-petrova-wojtaszczyk-2010}.)
\label{lem:cohen} 
If for some $I_C\in\mathbb{N}$ there exists $0\le I\le I_C$, $q,m\in\mathbb{N}$ and $\theta\in(0,1)$ 
such that $\theta \sigma_I \le \sigma_{I+qm}$ then it also holds $\sigma_I \le d_{m}/|\theta-q^{-\frac12}| $.
Moreover,
if %for some $I\in\mathbb{N}$ there exists $0\le i\le I$, $q,m\in\mathbb{N}$ and $\theta\in(0,1)$ 
%such that $\theta \tilde\sigma_i \le \tilde \sigma_{i+qm}$ 
the same hypothesis holds for $\theta \tilde\sigma_I \le \tilde \sigma_{I+qm}$,
then it holds $\tilde \sigma_I \le (d_{m}+\theta_{I_C}d_{I_C})/|\theta-q^{-\frac12}|$ with a probability more than $1-\epsilon$, more precisely in the events where it also holds $\tilde\sigma_I-\theta_Id_I \le V(Z^{\tilde\lambda_{I+1}}-\hat Y^{\tilde\lambda_{I+1}}) \le \tilde\sigma_I$.
\end{lemma}

Let us treat the case in the weak greedy algorithm where for $\alpha>0$ it holds $d_I\le d_0 I^{-\alpha}\,,\ \forall 0\le I\le I_M$. Assume we use the MC estimator $M_{\rm test}^{\epsilon,I_M}$ whatever $\epsilon\in(0,1)$.
For $I\le I_M$, the outcome $(\tilde\sigma)_{0\le i\le I_M}$ of the weak greedy algorithm is as follows
with probability more than $1-\epsilon$.
For $q\in\N_{>0}$, and $\beta,\gamma\in\R_{>0}$ to be determined later, we define $I_0:=\lceil\gamma q\rceil>0$
and assume
\begin{equation}
\label{eq:ababsurdo}
\forall C>0, \exists I_0\le \hat I\le I_M \slash\tilde \sigma_{\hat I}>C d_0 \hat I^{-\alpha} \,.
\end{equation}
Given $C>0$, we denote the smallest $\hat I$ satisfying~\ref{eq:ababsurdo} by $I_C$.
If for $0< I\le I_C$ we can define $m_I\in\N$ by $I+qm_I\le I_C\le I+q(m_I+1)$ then, since $(\sigma_i)_i$ is a decreasing sequence, it holds by~\eqref{eq:ababsurdo},
\begin{equation}
\label{ass:lemma}
\tilde \sigma_I \le \tilde \sigma_{I_C} (I_C/I)^{\alpha} \le \tilde \sigma_{I_{I+qm_I}} (I_C/I)^{\alpha} \,.
\end{equation}
Applying Lemma~\ref{lem:cohen} with $m=m_I$ and $\theta=(I/I_C)^{\alpha}$
\begin{equation}
%\begin{multline}
\tilde \sigma_{I_C} \le \tilde \sigma_{I} \le d_{m_I} (1+\theta_{I_C})/|(I/I_C)^{\alpha}-q^{-\frac12}| 
%\end{multline}
\end{equation}
(observe $d_{m_I}\ge d_{I_C}$) contradicts~\eqref{eq:ababsurdo} as soon as one chooses
$C\ge (I_C/m_I)^{\alpha}(1+\theta_{I_C})/|(I/I_C)^{\alpha}-q^{-\frac12}|>0$. In particular, if we set $I=\Theta I_C$ with some well-chosen $\Theta\in(0,1-1/\beta-1/\gamma)$, then
$$
I_C \ge I_0 \ge \gamma q \ge \frac{q}{1-\Theta-\frac1\beta}
$$
and $\beta q m_I \ge\beta(I_C-I-q)\ge\beta((1-\Theta)I_C-q)\ge I_C $.
In fact, it is always possible to find $I\in\N_{>0}$ with $m_I\in\N_{>0}$ ($m_I\ge\lfloor\gamma/\beta\rfloor>0$ in particular) provided one takes $\gamma>1$ large enough and sufficiently smaller $\beta>1$. This requires
$C \ge (\beta q)^{\alpha}(1+\theta_{I_C})/|\Theta^{\alpha}-q^{-\frac12}| \,.$
And finally, it also holds $\forall 0\le I< I_0, \tilde \sigma_{I}< C d_0 I^{-\alpha}$ provided one chooses $C$ such that
$$
C\ge\max\{(\beta q)^{\alpha}(1+\theta_{I_C})/|\Theta^{\alpha}-q^{-\frac12}|,(\tilde\sigma_0/d_0)I_0^{\alpha}\}>0 \,,
$$
which can still be optimized as a function of $q>0$.
So the outcome of the weak greedy algorithm satisfies $\forall 0\le I< I_M, \tilde \sigma_{I}< C d_0 I^{-\alpha}$
with probability more than $1-\epsilon$.
A similar reasoning can be adapted from~\cite{binev-cohen-dahmen-devore-petrova-wojtaszczyk-2010} for the other cases.

 \section*{Acknowledgement} 
 Part of this work was completed while the author was an academic host at MATHICSE-ASN Chair, EPFL (Lausanne, Switzerland).
 Thanks to Professors J. Rappaz and M. Picasso for their hospitality.
% the organizers of / participants to the ICMS ``Uncertainty Quantification''  workshop (Edinburgh, 2010 may)
% for giving him nice opportunities of discussion, 
% and Professors A.T. Patera and T. Leli\`evre for their numerous encouragements and advice.
%We would also like to thank the referee for their careful reading of the paper. 

%\begin{acknowledgements}
%\end{acknowledgements}

%\bibliographystyle{ws-m3as}
%\bibliographystyle{amsplain}
%\bibliographystyle{plain}
%\bibliographystyle{plain}
%\bibliography{bayesian.bbl}
%\bibliography{myrefs}
%\bibliography{/home/boyaval/Desktop/Z/database/myrefs.bib}%EDF
%\bibliography{/home/sebastien/Bureau/Z/database/myrefs}%portable
%\bibliography{$HOME/Bureau/Z/database/myrefs}

\end{document}